\documentclass[twoside,11pt]{article}

\usepackage{fancyheadings}%{fancyhdr}
\usepackage{latexsym}
\usepackage{lastpage,engord,graphics}
\usepackage[dvips]{color}
\usepackage{epic}
\usepackage{eepic}
\usepackage{rotating}
\usepackage{amssymb}
\usepackage{fancybox}

\makeatletter
\let\savedoddhead=\@oddhead
\def\@oddhead{\hfil {\footnotesize MAPPINGS WITH MAXIMAL RANK}  \hfil\thepage }
\makeatother

\makeatletter
\let\savedevenhead=\@evenhead
\def\@evenhead{\thepage\hfil {\footnotesize C. ABREU-SUZUKI}  \hfil}
\makeatother

\title{Maximal Rank Maps Between Riemannian Manifolds With 
Bounded Geometry}
\author{C. Abreu-Suzuki}
\newcommand{\noin} {\ensuremath{\in\kern-0.77em|}\hspace{0.06in}}
\newcommand{\real} {\ensuremath{I\kern-0.37emR}}
\newcommand{\Na} {\ensuremath{I{\kern-0.37emN}}}
\newcommand{\Z} {\ensuremath{Z\kern-0.40emZ}} % aqui!!!
\newcommand{\rn} {\ensuremath{I\kern-0.37emR^{n}}}
\newcommand{\ra}[1] {\ensuremath{I\kern-0.37emR^{{#1}}}}
\newcommand{\rnz} {\ensuremath{{I\kern-0.37emR}^{n}\setminus\{0\}}}

\newcommand{\pfe} {\mbox{ } \hfill$\Box$}

\newcommand{\pf} {\mbox{\it Proof.\ }}

\newcommand{\cl}[1] {{\it (Claim {#1})}}

\newtheorem{def21}{Definition}[section]
\newtheorem{def22}[def21]{Definition}
\newtheorem{prop23}[def21]{Proposition}
\newtheorem{def24}[def21]{Definition}

\newtheorem{def314}{Definition}[section]
\newtheorem{thm315}[def314]{Theorem}
\newtheorem{lem31}[def314]{Lemma}
\newtheorem{lem32}[def314]{Lemma}
\newtheorem{def313}[def314]{Definition}
\newtheorem{def33}[def314]{Definition}
\newtheorem{prop34}[def314]{Proposition}
\newtheorem{prop35}[def314]{Proposition}

\newtheorem{def37}[def314]{Definition}
\newtheorem{prop38}[def314]{Proposition}
\newtheorem{prop39}[def314]{Proposition}
\newtheorem{prop310}[def314]{Proposition}
\newtheorem{prop311}[def314]{Proposition}
\newtheorem{prop312}[def314]{Proposition}

\newtheorem{def54}{Definition}[section]
\newtheorem{def55}[def54]{Definition}
\newtheorem{def56}[def54]{Definition}
\newtheorem{lem51}[def54]{Lemma}
\newtheorem{prop52}[def54]{Proposition}
\newtheorem{thm53}[def54]{Theorem}

\makeindex
\begin{document}

\date{\empty}

%\maketitle
\thispagestyle{plain}

\begin{center}
\textbf{\large MAXIMAL RANK MAPS BETWEEN RIEMANNIAN MANIFOLDS WITH 
BOUNDED GEOMETRY}\footnote{\textit{Mathematics Subject Classification.} Primary 53C20.}
\end{center}

\begin{center}
{C. ABREU-SUZUKI}
\end{center}

%\begin{abstract}
%If $\pi: M\rightarrow B$ is an onto smooth maximal rank map  
%between complete Riemannian manifolds $M$ and $B$ with bounded 
%geometry, we prove sufficient conditions for $M$ to be 
%roughly isometric to the Riemannian product $F\times B$, where 
%$F$ is a fiber of $M$.
%\end{abstract}

   \pagenumbering{arabic}

\section{Introduction}
\protect\label{ls1}
\thispagestyle{plain}

Rough isometries, in the sense of M. Kanai~\cite{MK1}, 
provide equivalence relations between non-compact Riemannian manifolds.
M. Kanai showed that when two spaces are roughly isometric they share 
properties such  as volume growth rate and the validity of 
isoperimetric inequalities. He accomplished that via approximating 
a Riemannian manifold by a combinatorial structure, he calls a net. 
He proved that complete Riemannian manifolds, whose Ricci curvature 
are bounded from below, are roughly isometric to nets. We provide 
background in section~\ref{ls2}. 

Here we study mappings with maximal rank $\pi: M\rightarrow B$, 
between complete non-compact Riemannian manifolds $M$ and $N$ with 
bounded geometry.  
O'Neill~\cite{BON} gives necessary and sufficient conditions for a 
Riemannian submersion $\pi: M\rightarrow B$ to be trivial, 
i.e., to differ only by an isometry of $M$ from the simplest type of 
Riemannian submersions, the projection $p_{B}: F\times B \rightarrow B$
 of a Riemannian product manifold $F\times B$ on one of its factors $B$ 
(see \textbf{Theorem~\ref{lthm315}}). 
In section~\ref{ls3} we review O'Neill's 
results and describe the properties of long curves in $B$ lifted to 
$M$.

In section~\ref{ls5} we define two new properties of maximal rank onto 
mappings $\pi: M\rightarrow B$: uniformly 
roughly isometric fibers [{\bf Definition~\ref{ldef55}}] and 
horizontal lift control [{\bf Definition~\ref{ldef56}}].
Then we prove that if $M$ and $B$ are complete Riemannian manifolds  
with bounded geometry, and if $\pi$ satisfies these two properties
with trivial holonomy, then $M$ is roughly isometric to the 
product $F\times B$ of the base manifold $B$ and a fixed fiber $F$ of 
$M$ [{\bf Theorem~\ref{lthm53}}].

\section{Rough Isometries, Nets and Bounded Geometry}
\protect\label{ls2}

In this section we introduce notation, give a few definitions
according to M.Kanai~\cite{MK1}  
and O'Neill~\cite{BON}, and state some
results without proofs, providing  references whenever necessary.

Rough isometries, a concept first introduced by M.  Kanai~\cite{MK1} give 
equivalence relations, which will be of our interest.

\begin{def21} 
\protect\label{ldef21}
Let $(M,\delta)$ and $(N,d)$ be metric spaces.
A map $\varphi : M\rightarrow~N$, not necessarily continuous, is
called a \textbf{rough isometry}, if it satisfies the 
following two axioms:
\begin{description}
\item[(RI.1)] There exist constants $A\geq 1, C\geq 0$, satisfying, 
\\
\[
\frac{1}{A}\delta(p_{1},p_{2}) - C \leq
d(\varphi(p_{1}),\varphi(p_{2}))
\leq A
\delta(p_{1},p_{2}) + C, \hspace{0.2in} \forall p_{1},p_{2}\in M 
\]
\item[(RI.2)] The set $Im\varphi := \{ q = \varphi(p), \forall p \in
M\}$ is
\textbf{full} in $N$, i.e.  
\[
\exists\varepsilon >0: N = B_{\varepsilon}(Im\varphi) =
\{q \in N: d(q,Im\varphi)<\varepsilon\}
\] 
In this case we say that $Im\varphi$ is \textbf{$\varepsilon$-full} 
in $N$.
\end{description}
\end{def21}

It is immediate to verify that if $\varphi : M\rightarrow N$ 
and $\psi : N\rightarrow M$ are rough isometries, then  
the composition $\varphi\circ\psi : N\rightarrow N$ is also a 
rough isometry.

A rough inverse of $\varphi$, which we will denote by 
$\varphi^{-} : N \rightarrow M$ is defined as follows: 
for each $q\in N$, choose $p\in M$ such that 
$d(\varphi(p),q)<\varepsilon$. Such a $p$ exists because of axiom  
{\bf (RI.2)}.  $\varphi^{-}$ is a rough isometry such that both 
$\delta(\varphi^{-}\circ\varphi(p),p)$ and 
$d(\varphi\circ\varphi^{-}(q),q)$
are bounded in $p\in M$ and in $q\in N$, respectively.

To study geometric properties of manifolds, which are invariant under
rough isometries, we next introduce what is called  in~\cite{MK1}, a 
\textbf{net}. A net  is a discrete or combinatorial structure 
that provides approximations of Riemannian manifolds.

\begin{def22} 
\protect\label{ldef22}
 Let $P$ be a countable set. A family $N=\{N(p):p\in P\}$ is
called a \textbf{net structure} of $P$ if the following conditions
hold for all $p,q\in P$:
    \begin{description}
    \item[(N.1)] $N(p)$ is a finite subset of $P$
    \item[(N.2)] $q\in N(p)$ iff $p\in N(q)$
    \end{description}
\end{def22}

Let $M$ be a complete Riemannian manifold, and let $d$ be the
induced metric. A subset $P$ of $M$ is said to be
\textbf{$\varepsilon$-separated} for $\varepsilon >0$, 
if $d(p,q)\geq\varepsilon$ whenever $p$ and $q$ are distinct 
points of $P$, and an
$\varepsilon$-separated set is called \textbf{maximal} if it 
is maximal with respect to the order relation of inclusion. 

We have the following, 

\begin{prop23}
\protect\label{lprop23}
If $P$ is a countable maximal $\varepsilon-$separated set in a 
Riemannian manifold $(M,d)$, then $P$ is $\varepsilon-$full in $M$, 
where $\varepsilon >0$.
\end{prop23}
\pf\hspace{0.1in}
We want to show that,
\[
d(x,P)<\varepsilon, \hspace{0.2in}\forall x\in M
\]

If $x\in P$, then $d(x,P)=0<\varepsilon$.

If $x\in M\setminus P$, by the maximality of $P$, there exists 
$\bar{p}\in P$ such that $d(x,\bar{p})<\varepsilon$, and finally the 
definition of infimum implies that 
$d(x,P):=\inf_{p\in P}d(x,p)\leq d(x,\bar{p})<\varepsilon.$ 

\pfe

Let $P$ be a maximal $\varepsilon$-separated subset of $M$. 
We define a \textbf{net structure} $N=\{N(p):p\in P\}$ of $P$ by 
$N(p)=\{q\in P:0<d(p,q)\leq 2\varepsilon\}$. A maximal
$\varepsilon$-separated subset of a  complete Riemannian 
manifold with the net structure described above will be called an
\textbf{$\varepsilon$-net in $M$}.

For a point $p\in P$, each element of $N(p)$ is called a 
\textbf{neighbor}
of $p$. A sequence $p=(p_{0},\cdots,p_{l})$ of points in $P$ is called 
a {\it path} from $p_{0}$ to $p_{l}$ of length $l$ if each $p_{k}$ is a 
neighbor of $p_{k-1}$. A net $P$ is said to be \textbf{connected} 
if any 
two points in $P$ are joined by a path. For points $p$ and $q$ of a 
connected net $P$, $\delta(p,q)$ denotes the minimum of the lengths of
paths from $p$ to $q$. This $\delta$ satisfies the axioms of metric and 
it is called, according to~\cite{MK1}, the 
\textbf{combinatorial metric} of $P$.

We observe that an $\epsilon-$net in a complete Riemannian manifold is 
connected if the manifold is connected (see~\cite{MK1}). 

In what follows, we introduce some notation
(c.f.~\cite{MK1})
and we define a \textbf{bounded geometry} condition
for manifolds. 

\vspace{0.1in}

Let $(M,g)$ be a Riemannian manifold with Levi-Civita connection
$\nabla$ and curvature tensor $R$.

The \textbf{Ricci curvature tensor} of $(M,g)$, at each $x \in M$ 
is a symmetric bilinear form $Ric$ defined by
\[
\begin{array}{rcll}
Ric: & T_{x}M\times T_{x}M & \longrightarrow &\real \\
     & (\xi,\mu)           & \longmapsto     &
     Ric(\xi,\mu):=trace(\zeta\mapsto R(\xi,\zeta)\mu)
\end{array}
\]
If $M$ is complete, the \textbf{injectivity radius at} $x\in M$ 
is given by
\[
\imath_{x}(M):=\sup\{r>0: \exp_{x}\mid_{B(x,r)}
\mbox{ is a diffeomorphism}\}
\]
and $\imath(M):=\inf\{\imath_{x}(M): x\in M\}$
is called the {\it injectivity radius} of $M$.

\begin{def24}
\protect\label{ldef24}
Let $M$ be a complete $m$-dimensional Riemannian manifold.
We say that $M$ has \textbf{bounded geometry} if it satisfies:
\begin{description}
\item[(BG.R)] the Ricci curvature is bounded from below by
$-(m-1)k_{M}^{2}$,  where $k_{M}$ is a positive constant;
\item[(BG.I)] the injectivity radius $\imath(M)$ is positive.
\end{description}
\end{def24}

We recall that a complete Riemannian manifold satisfying a 
bounded geometry condition has its geometry reflected by that of any 
net that approximates the manifold (see~\cite{MK1}, Lemma2.5).

\section{Long Curves and O'Neill Diffeomorphisms}
\protect\label{ls3}

Here we review background from O'Neill~\cite{BON} and 
Abreu-Suzuki~\cite{CASMMR} concerning mappings of maximal rank.

Let  $M$ and $B$ be Riemannian manifolds with dimensions $m$ and $n$,
respectively, where $m\geq n$.
We will denote by $\pi: M\rightarrow B$ an onto mapping with maximal rank $n$, 
that is, $\pi$ and each of its derivative maps $\pi_{\ast}$  are surjective.  

We start recalling the definitions of horizontal and vertical vectors, and  
of a Riemannian submersion, according to~\cite{BON}. 

A  tangent vector on $M$ which is tangent to a fiber is called 
\textbf{vertical}, and if it is orthogonal to  a fiber it is called 
\textbf{horizontal}. So, if a vector field  on $M$ is always tangent to  
fibers, we say that it is {vertical}, and  if it is always orthogonal to  
fibers, we say that it is {horizontal}. 

\begin{def314}
\protect\label{ldef314}
A \textbf{Riemannian} \textbf{submersion} 
$\pi: M\rightarrow B$ is an onto mapping, such that, $\pi$ has maximal 
rank, and  $\pi_{\ast}$ preserves lengths of horizontal vectors.
\end{def314}

Because for all $x\in M$ each derivative map $\pi_{\ast}(x)$ of 
$\pi$ is surjective, we can define the projections ${\cal H}$ and ${\cal V}$ 
of the tangent space of $M$ onto the subspaces of horizontal and vertical 
vectors, respectively, which will be denoted, respectively by $(VT)_{x}$ 
and $(HT)_{x}$ for each $x\in M$.
In that case, we can decompose each tangent space to $M$ into the direct
orthogonal sum $T_{x}M=(VT)_{x}\oplus (HT)_{x}$.

Recall, O'Neill\cite{BON} proved,

\begin{thm315}
\protect\label{lthm315} \textbf{(O'Neill)}  
Let $\pi: M\rightarrow B$ 
be a submersion of a complete Riemannian manifold $M$. Then $\pi$ is 
trivial if and only if the tensor $T$  and the 
group $G$ of the submersion both vanish. 
\end{thm315}

O'Neill defines the tensor $T$ on $M$, which is the second fundamental form 
of all fibers, by $T_{E}F = {\cal H}\nabla _{{\cal V}E}({\cal V}F) + 
{\cal V}\nabla _{{\cal V}E}({\cal H}F)$ for arbitrary vector fields $E$ 
and $F$, where $\nabla$ is the covariant derivative of $M$. The group $G$ 
of the submersion is the holonomy group of the connection 
$\Gamma (x\in M\mapsto {\cal H}(T_{x}M)=(HT)_{x})$, with reference to the 
base point $O\in B$.

The unique horizontal vector property, as stated in 
{\bf Lemma~\ref{llem31}}, follows from the maximality of the rank of 
the onto mapping $\pi$.

\begin{lem31}
\protect\label{llem31}
 Let $b\in B$ be fixed. For any $w\in T_{b}B$ and $x\in M$
such that $\pi(x)=b$, there exists  a unique horizontal vector
 $v\in T_{x}M$ which is $\pi$-related to $w$, i.e. satisfying  $v\in
(HT)_{x}$ and  $(\pi_{\ast})_{x}(v)=w$.%doesn't use S.2
\end{lem31}

In the following Lema, with additional control from below over 
the length  of horizontal vectors, one has control from below 
over the distance in $M$. 

\begin{lem32}
\protect\label{llem32}
Let $M$ and $B$ be connected and geodesically complete.
For any $x,x'\in M$, let $\Gamma_{min}\subset M$ be a minimal geodesic 
joining $x$ to $x'$, and let $\gamma_{min}\subset B$ be a minimal 
geodesic joining $\pi(x)$ to $\pi(x')$.
Assume that for all $b \in B$ and for all
$x\in F_{b}$ there exist constants $\alpha\geq 1$ and $\beta >0$,
both independent of $b$ and $x$, satisfying
\begin{equation}
\protect\label{lequ47}
\frac{1}{\alpha} ||w||_{B} - \beta \leq ||v||_{M}
\end{equation}
for all $w\in T_{b}B$,  where $v$ is the unique horizontal lift of
$w$ through $x$ that we assume satisfies $||v||_{M}\leq 1$, where 
$||\hspace{0.15in}||_{M}$, $||\hspace{0.15in}||_{B}$ denote the
inner product on $TM$ and $TB$, respectively.

Then, \hspace{0.1in}
$
d_{M}(x,x')=
\ell(\Gamma_{min})\geq \frac{1}{\alpha} \ell(\gamma_{min})-\beta
=\frac{1}{\alpha}d_{B}(\pi(x),\pi(x'))-\beta
$
\end{lem32}
\pf\hspace{0.1in}
The proof of this Lemma is in~\cite{CASMMR} .
 
\vspace{0.1in}

Recall the definition of a lift of a curve.

\begin{def313}
\protect\label{ldef313}
Let $\gamma:[t_{1},t_{2}]\rightarrow B$ be a smooth embedded curve 
in $B$. A curve $\Gamma:[t_{1},t_{2}]\rightarrow M$ 
satisfying $\pi\circ\Gamma=\gamma$ is called a \textbf{lift} of $\gamma$.

If in addition, $\Gamma$ is horizontal, 
i.e., ${\Gamma}^{\prime}(t)\in (HT)_{\Gamma(t)}, 
\forall t\in [t_{1},t_{2}]$, where 
$\Gamma(t_{1})=x_{0}\in M$ with $\gamma(t_{1})=\pi(x_{0})$, the curve 
$\Gamma$ is called a \textbf{horizontal lift} of $\gamma$ 
through $x_{0}$.
Recall that the horizontal lift of a curve in $B$, through a 
point $x_{0}\in M$ is unique.
\end{def313}

\vspace{0.1in}

We now define long curves~\cite{CASMMR}.

\begin{def33}
\protect\label{ldef33}
Let $\beta>0$ be any fixed constant. A smooth embedded curve 
$\gamma:[t_{1},t_{2}]\rightarrow B$ is said to be a 
\textit{$\beta$-long curve} if 
$\inf_{t_{1}\leq t\leq t_{2}}||\gamma'(t)||\geq\beta$. In that case, 
$\ell(\gamma)\geq \displaystyle{\int}_{t_{1}}^{t_{2}}||\gamma'(t)|| dt 
\geq \beta (t_{1}-t_{2})$. We say that a curve $\gamma$ is simply a 
long curve if it is a $\beta$-long curve for some constant $\beta >0$.
\end{def33}

Let $\gamma:[t_{1},t_{2}]\rightarrow B$ denote a smooth embedded curve
and let $\Gamma:[t_{1},t_{2}]\rightarrow M$ denote a lift of $\gamma$.

In the next two  Propositions, proven in~\cite{CASMMR}, under control 
from above (or below)  on the derivative of the maximal rank mapping 
$\pi$, we have control from below (or above) over the length of any lift 
of a curve. 
In {\bf Proposition~\ref{lprop34}} any lift $\Gamma$ in $M$
of a long curve $\gamma$ in $B$  cannot be short, and  in 
{\bf Proposition~\ref{lprop35}}
the length of a lift $\Gamma$ of a long curve $\gamma$ is
bounded above by the  length of $\gamma$.

The Riemannian norms in $TM$ and $TB$ will be denoted by 
$||\hspace{0.1in}||_{M}$ and $||\hspace{0.1in}||_{B}$, respectively.

\begin{prop34}
\protect\label{lprop34}
Let $\alpha\geq 1$ and $\beta >0$ be constants satisfying,
\begin{equation}
\protect\label{lequ42}
\left||(\pi_{\ast})_{x}v|\right|_{B} \leq
\alpha\left||v|\right|_{M}+\beta
\end{equation}
for all $x\in M$, for all $v \in T_{x}M$ where 
$\left||v|\right|_{M}\leq 1$.

If $\gamma$ is any smooth $\beta$-long curve in $B$, then,
\[
\ell(\Gamma)\geq\frac{1}{\alpha}\left[\ell(
\gamma)-\beta (t_{2}-t_{1})\right]>0
\]
where $\ell(\Gamma)$ and $\ell(\gamma)$ denote the lengths of
the curves $\Gamma$ and $\gamma$, respectively.
\end{prop34}

\begin{prop35}
\protect\label{lprop35} Let $\Gamma$ be a lift of $\gamma$, and  
assume for \underline{horizontal}$^{(\dagger)}$
%\underline{non-vertical}$^{(\dagger)}$
vectors $v\in TM$ only,
that there is a universal constant $\alpha\geq 1$
%are universal constants $\alpha\geq 1, \beta >0$
satisfying,
\begin{equation}
\protect\label{lequ46}
\left||(\pi_{\ast})_{x}v|\right|_{B} \geq
\frac{1}{\alpha}\left||v|\right|_{M}-\beta
\end{equation}
for all $x\in M$, for all $v \in T_{x}M\setminus (VT)_{x}
=(HT)_{x}=\left[\ker (\pi_{\ast})_{x}\right]^{\perp}$.

For a $\beta$-long curve $\gamma$, we have,
\[
\ell(\Gamma)\leq \alpha\left[\ell(\gamma)+\beta (t_{2}-t_{1})\right]
\]
where $\ell(\Gamma)$ and $\ell(\gamma)$ denote the lengths of
the curves $\Gamma$ and $\gamma$, respectively.
\end{prop35}

Next, for onto smooth mappings with maximal rank between complete 
Riemannian manifolds,  we recall the definition of 
special diffeomorphisms between any two fibers, we call 
\textbf{O'Neill diffeomorphisms}, a useful tool that will feature 
in many of our proofs. %{\bf Theorem~\ref{lthm412}} 

Firstly, we give the definition, and in the five propositions 
that follow we state several of their properties
 (see~\cite{CAS} Theorem 4.12).

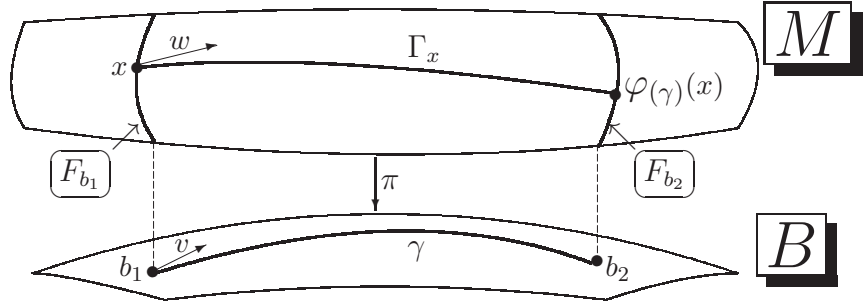
\begin{figure}[here]

         %\begin{picture}(100,360) (0,0)
            %\includegraphics{picthm412-1.ps}
         \begin{picture}(370,110)(0,0)

%\dottedline{2}(0,120)(370,120)
%\dottedline{2}(0,0)(370,0)

\put(330,90){\shadowbox{\Huge $M$}}
\qbezier(50,110)(180,120)(320,110)
\qbezier(50,70)(200,50)(320,70)
\qbezier(50,70)(40,85)(50,110)
\qbezier(320,70)(335,95)(320,110)

\put(60,50){\ovalbox{\large $F_{b_{1}}$}}
\put(83,64){$\nearrow$}
\thicklines
\qbezier(99,65)(85,85)(99,113)
\thinlines
\put(90,90){$\bullet$}
\put(83,90){$x$}
\put(93,94){\vector(4,1){30}}
\put(105,100){$w$}

\put(281,50){\ovalbox{\large $F_{b_{2}}$}}
\put(271,64){$\nwarrow$}
\thicklines
\qbezier(267,63)(282,95)(267,113)
\thinlines
\put(271,80){$\bullet$}
\put(277,83){{\Large $\varphi_{(\gamma)}$}$(x)$}

\thicklines
\qbezier(93,93)(160,100)(274,83)
\thinlines
\put(195,97){\large $\Gamma_{x}$}

\put(327,10){\shadowbox{\Huge $B$}}
\qbezier(53,15)(186,60)(320,15)
\qbezier(103,5)(186,15)(270,5)
\qbezier(53,15)(65,15)(103,5)
\qbezier(270,5)(300,15)(320,15)

\thicklines
\qbezier(99,15)(195,45)(264,19)
\thinlines
\put(195,22){\large $\gamma$}
%\put(205,23){\large $\leadsto$}
\put(96,13){$\bullet$}
\put(86,14){$b_{1}$}
\dashline[+90]{3}(99,17)(99,65)
\put(99,16){\vector(2,1){20}}
\put(107,23){$v$}
\put(264,17){$\bullet$}
\put(270,14){$b_{2}$}
\dashline[+90]{3}(267,19)(267,64)

%text
\put(183,59){\vector(0,-1){20}}
       \put(185,47){\large $\pi$}

         \end{picture}

      \caption{The map $\varphi_{(\gamma)}$.}
         \label{fig1def37}
         \index{pictures!Theorem\ref{lprop37}}
      \end{figure}

\begin{def37}
\protect\label{ldef37}
Let $\pi: M\rightarrow B$ be an onto smooth map with
maximal rank, where $M,B$ are complete and $B$ is connected.
Let  $b_{1}, b_{2}$  be distinct elements of $B$ and
let $\gamma:[t_{1},t_{2}]\longrightarrow B$ be a piecewise smooth
embedded curve parametrized proportionally to arclength, 
where $\gamma(t_{1})=b_{1}, \gamma(t_{2})=b_{2}$.

If we denote by $\pi^{-1}(b_{1}) = F_{b_{1}}$ and
$\pi^{-1}(b_{2}) = F_{b_{2}}$ their corresponding fibers, we thus
define the map  
$ \varphi_{(\gamma)}:F_{b_{1}}\longrightarrow  F_{b_{2}}$, 
we refer to as \textbf{O'Neill diffeomorphism},
by the following rule:
Given $x \in F_{b_{1}}$, let $\Gamma_{x}$
be the unique horizontal lift of $\gamma$ through $x$, and set
$\varphi_{(\gamma)}(x):= \Gamma_{x}(t_{2}) \in F_{b_{2}}.$
(see Fig.~{\ref{fig1def37}})
\end{def37}

%(see Fig.~{\ref{fig1def37}})

O'Neill~\cite{BON} noted that his diffeomorphisms have the following 
five properties 
(detailed proofs are available in author's thesis~\cite{CAS}):

\begin{prop38}
\protect\label{lprop38}
 $\varphi_{(\gamma)}: F_{b_{1}}\longrightarrow F_{b_{2}}$ is
 well-defined.
\end{prop38}

\begin{figure}[here]

         \begin{picture}(370,95)(0,0)

%\dottedline{2}(0,120)(370,120)
%\dottedline{2}(0,0)(370,0)

\put(330,90){\shadowbox{\Huge $M$}}
\qbezier(50,110)(180,120)(320,110)
\qbezier(50,70)(200,50)(320,70)
\qbezier(50,70)(40,85)(50,110)
\qbezier(320,70)(335,95)(320,110)

\thicklines
\qbezier(99,65)(85,85)(99,113)
\thinlines

\put(97,94){\vector(1,0){55}}
\put(110,102){{\large $\varphi_{(\gamma_{1})}$}}

\thicklines
\qbezier(160,61)(150,85)(160,115)
\thinlines

\put(160,94){\vector(1,0){113}}
\put(200,102){{\large $\varphi_{(\gamma_{2})}$}}

\thicklines
\qbezier(267,63)(282,95)(267,113)
\thinlines

\put(99,84){\vector(1,0){170}}
\put(160,75){{\large $\varphi_{(\gamma_{2}\circ\gamma_{1})}$}}

\put(325,10){\shadowbox{\Huge $B$}}
\qbezier(53,15)(186,60)(320,15)
\qbezier(103,5)(186,15)(270,5)
\qbezier(53,15)(65,15)(103,5)
\qbezier(270,5)(300,15)(320,15)

\put(96,12){$\bullet$}
\dashline[+90]{3}(99,17)(99,65)
\thicklines
\qbezier(99,15)(137,35)(158,13)
\thinlines
\put(126,15){\large $\gamma_{1}$}
%\put(205,23){\large $\leadsto$}
\put(157,11){$\bullet$}
\dashline[+90]{3}(160,17)(160,61)
\thicklines
\qbezier(158,13)(200,40)(264,19)
\thinlines
\put(205,20){\large $\gamma_{2}$}
\put(264,16){$\bullet$}
\dashline[+90]{3}(267,19)(267,64)

%text
\put(183,59){\vector(0,-1){20}}
       \put(185,47){\large $\pi$}

         \end{picture}

      \caption{Property $\varphi(\gamma_{2}\circ\gamma_{1})=
      \varphi(\gamma_{2})\circ\varphi(\gamma_{1})$.}
         \label{fig2prop39}
         \index{pictures!Theorem\ref{lprop39}}
      \end{figure}
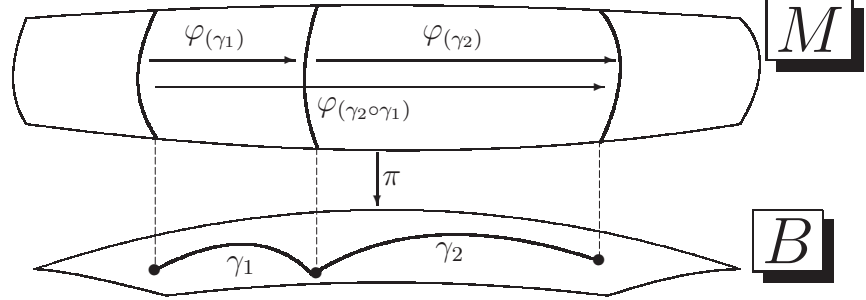

\begin{prop39}
\protect\label{lprop39}
 Let $\gamma_{1}:[t_{1},t_{2}]\longrightarrow B$ and
$\gamma_{2}:[t_{2},t_{3}]\longrightarrow B$ be smooth embedded curves
parametrized proportionally to arclength, such that 
$\gamma_{1}(t_{2}) = \gamma_{2}(t_{2})$.
Define $\gamma_{3}: [t_{1},t_{3}] \longrightarrow  B$  the composition 
of $\gamma_{1}$  and $\gamma_{2}$,
denoted by $\gamma_{3} = \gamma_{2}\circ\gamma_{1}$,
as follows: 
$
 t\in [t_{1},t_{3}]  \longmapsto  \gamma_{3}(t):=
     \left\{
     \begin{array}{ll}
       \gamma_{1}(t), & \mbox{ if } t_{1}\leq t\leq t_{2}\\
       \gamma_{2}(t), & \mbox{ if } t_{2}\leq t\leq t_{3}
     \end{array}
     \right.%\}
$ %\\
Then $\varphi_{(\gamma_{3})}:F_{\gamma_{1}(t_{1})}\longrightarrow
F_{\gamma_{2}(t_{3})}$ satisfies $\varphi_{(\gamma_{3})} =
\varphi_{(\gamma_{2})}\circ\varphi_{(\gamma_{1})}$. 
(see Fig.~{\ref{fig2prop39}})
\end{prop39}

%(see Fig.~{\ref{fig2prop39}})

\begin{prop310}
\protect\label{lprop310}
 $\varphi_{(\gamma)}: F_{b_{1}}\longrightarrow F_{b_{2}}$  is a
 diffeomorphism.
\end{prop310}

\begin{figure}[here]

         \begin{picture}(370,185)(0,0)

%\dottedline{2}(0,210)(370,210)
%\dottedline{2}(0,0)(370,0)

\put(330,170){\shadowbox{\Huge $M$}}
\qbezier(50,190)(180,205)(320,190)
\qbezier(50,110)(200,90)(320,110)
\qbezier(50,110)(30,150)(50,190)
\qbezier(320,110)(335,150)(320,190)

\put(82,89){\ovalbox{\large $F_{b_{1}}$}}
\put(105,103){$\nearrow$}
\thicklines
\qbezier(119,103)(105,145)(119,195.5)
\thinlines
\put(110,155){$\bullet$}
\put(103,158){$x$}

\put(260,88){\ovalbox{\large $F_{b_{2}}$}}
\put(250,102){$\nwarrow$}
\thicklines
\qbezier(247,102)(262,145)(247,195.5)
\thinlines
\put(251,160){$\bullet$}
\put(268,160){\ovalbox{$\varphi_{(\gamma)}(x)$}}
\put(257,160){$\leftarrow$}
\put(252,143){$\bullet$}

\put(160,169){$\Gamma_{x}$}
\thicklines
\qbezier(113,158)(180,170)(254,163)
\thinlines

%lift of arbitrary curve
\put(220,155){$\Lambda_{x}$}
\qbezier(113,158)(183,165)(255,146)
\put(181,160){\oval(170,40)}%here!!!!!!old=(183,160)
\put(183,126.5){\vector(0,1){13}}
\put(150,115){\ovalbox{$\forall$ \mbox{\scriptsize tub. neigh.}
${\cal T}_{\Gamma_{x}}$}}

\put(327,20){\shadowbox{\Huge $B$}}
\qbezier(53,35)(186,80)(320,35)
\qbezier(103,5)(186,20)(270,5)
\qbezier(53,35)(65,35)(103,5)
\qbezier(270,5)(300,35)(320,35)

\put(160,39){$\gamma$}
\thicklines
\qbezier(119,26)(175,45)(244,30)
\thinlines
%\put(205,23){\large $\leadsto$}
\put(116,23){$\bullet$}
\put(106,24){$b_{1}$}
\dashline[+90]{3}(119,27)(119,103)
\put(244,27){$\bullet$}
\put(251,28){$b_{2}$}
\dashline[+90]{3}(247,29)(247,102)
%arbitrary curve
\put(220,24){$\lambda$}
\qbezier(119,25)(150,35)(180,20)
\qbezier(180,20)(212,15)(245,29)
\put(183,30){\oval(165,30)}
\put(183,61.5){\vector(0,-1){16}}
\put(150,68){\ovalbox{$\exists$ \mbox{\scriptsize tub. neigh.}
${\cal T}_{\gamma}$}}

%text
\put(303,105){\vector(0,-1){62}}
       \put(305,70){\Large $\pi$}

         \end{picture}

      \caption{$\varphi_{(\gamma)}$ depends continuously on $\gamma$.}
         \label{fig3prop311}
         \index{pictures!Proposition\ref{lprop311}}
      \end{figure}
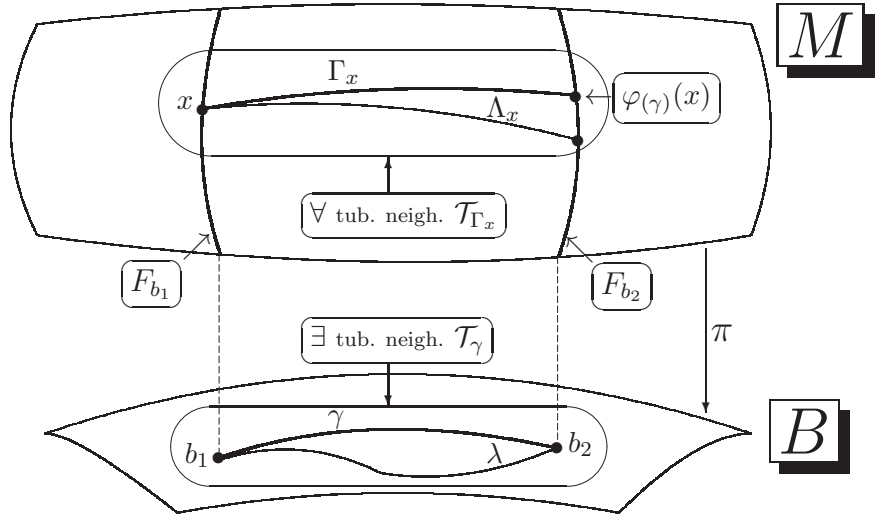

\begin{prop311}
\protect\label{lprop311}
$\varphi_{(\gamma)}$ depends continuously on $\gamma$, i.e.,
for sufficiently small displacements of $\gamma$ within a tubular
neighborhood, keeping the endpoints fixed, their horizontal lifts
through $x$ lie entirely within any given small tubular neighborhood of
$\Gamma_{x}$. 
Furthermore, since horizontal lifts are, by definition,
curves tangent to $(HT)\subset TM$, this fact forces the endpoints of
horizontal lifts through $x$ to belong to arbitrary balls around
$\varphi_{(\gamma)}(x)= \Gamma_{x}(t_{2})$, as long as those horizontal
lifts through $x$ lie entirely within a sufficiently small tubular
neighborhood of $\Gamma_{x}$. (see Fig.~{\ref{fig3prop311}})
\end{prop311}

%(see Fig.~{\ref{fig3prop311}})

%A group of diffeomorphisms between fibers is defined as follows.

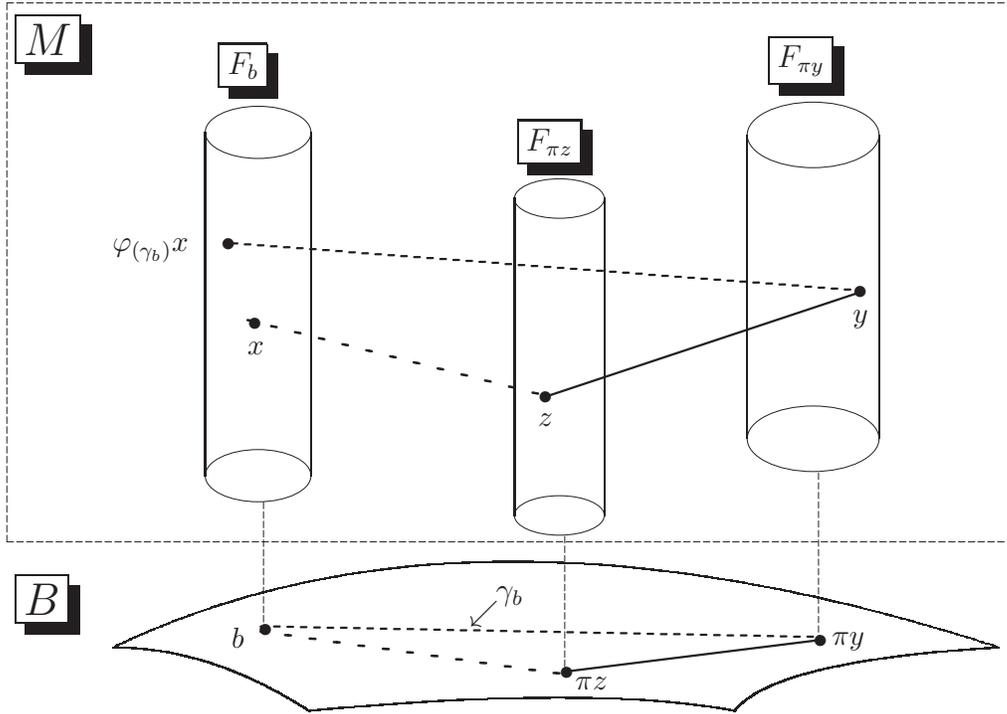
\begin{figure}[here]%oldthm

         \begin{picture}(390,270)(0,-15)%(390,270)(0,0)

%\graphpaper[20](-75,0)(470,275)%default spacing is [10]
\definecolor{darkyellow}{rgb}{0.7,0.7,0}
\definecolor{orange}{rgb}{1.0,0.4,0}
\definecolor{darkgreen}{rgb}{0,0.8,0}
%\dottedline{2}(0,300)(390,300)%(0,270)(390,270)
%\dottedline{2}(0,0)(390,0)

\put(8,242){\shadowbox{\LARGE $M$}}
\dashline[+90]{3}(5,269)(386,269)
\dashline[+90]{3}(5,65)(5,269)
\dashline[+90]{3}(386,65)(386,269)
\dashline[+90]{3}(5,65)(386,65)

\put(85,233){\shadowbox{\large $F_{b}$}}

%cylinder
%parallels=horizontal curves
\put(100,220){\ellipse{40}{20}}
\put(100,90){\ellipse{40}{20}}
%meridians=vertical curves
\put(80,90){\line(0,1){130}}
\put(120,90){\line(0,1){130}}
\dashline[+90]{3}(102,80)(102,31)

\put(198,204){\shadowbox{\large $F_{\pi z}$}}

%cylinder
%parallels=horizontal curves
\put(214,195){\ellipse{34}{15}}
\put(214,75){\ellipse{34}{15}}
%meridians=vertical curves
\put(197,75){\line(0,1){120}}
\put(231,75){\line(0,1){120}}
\dashline[+90]{3}(216,67)(216,14)

\put(293,234){\shadowbox{\large $F_{\pi y}$}}

%cylinder
%parallels=horizontal curves
\put(310,219){\ellipse{50}{25}}
\put(310,104){\ellipse{50}{25}}
%meridians=vertical curves
\put(285,104){\line(0,1){115}}
\put(335,104){\line(0,1){115}}
\dashline[+90]{3}(312,91)(312,27)

\put(8,30){\shadowbox{\LARGE $B$}}
\qbezier(45,25)(180,90)(380,25)
\qbezier(119,1)(227,11)(280,1)
\qbezier(45,25)(82,25)(119,1)
\qbezier(380,25)(295,20)(280,1)

\thicklines
\put(217,16){\line(8,1){97}}%line pi z-pi y
%\put(217,16){{\color{orange}\line(8,1){97}}}%line pi z-pi y
\dashline{3}[0.7](101,31)(214,15)%dashline b-pi z
%%\put(217,16){\line(-7,1){115}}%line b-pi z
%\put(217,16){{\color{green}\line(-7,1){115}}}%line b-pi z
\dashline[+90]{3}(101,33)(310,29)%dashline b-pi y
%%\put(104,32){\line(45,-1){208}}%line b-pi y
%\put(104,32){{\color{blue}\line(45,-1){208}}}%line b-pi y
\thinlines

\put(90,25){$b$}
%\put(102,21){\footnotesize $\in P_{B}$}
\put(100,29){$\bullet$}
\put(220,9){$\pi z$}
\put(214,13){$\bullet$}
\put(317,26){$\pi y$}
\put(310,25){$\bullet$}
\put(190,43){\large $\gamma_{b}$}
\put(180,35){$\swarrow$}

%text
\put(45,175){$\varphi_{({\gamma_{b}})}x$}
\put(86,175){$\bullet$}
\put(96,136){$x$}
\put(96,145){$\bullet$}
\thicklines
\dashline[+90]{3}(90,178)(329,160)
%{\color{blue}\dashline[+90]{3}(90,178)(329,160)}
\dashline{3}[0.7](96,149)(206,121)
%{\color{green}\dashline{3}[0.7](96,149)(206,121)}
\thinlines
\put(325,148){$y$}
\put(325,157){$\bullet$}
\put(206,109){$z$}
\put(206,117){$\bullet$}
\thicklines
\put(209,120){\line(3,1){120}}
%\put(209,120){{\color{orange}\line(3,1){120}}}
         \end{picture}

         \caption{Holonomy Group $G_{b}$.}
         \label{fig4prop312}
         \index{pictures!Proposition\ref{lprop312}}
      \end{figure}

\begin{prop312}
\protect\label{lprop312}
Let $b\in B$ be a fixed base point and 
$\gamma_{b}:[t_{1},t_{2}]\longrightarrow B$ be a piecewise smooth
embedded geodesic loop parametrized proportionally to arclength, 
where $\gamma_{b}(t_{1})=\gamma_{b}(t_{2})=:b$.  The properties of 
$\varphi_{(\gamma_{b})}: F_{b}\longrightarrow  F_{b}$, imply that 
the set of mappings $G_{b} :=\{\varphi_{(\gamma_{b})}: 
F_{b}\rightarrow  F_{b}, \hspace{0.1in}\forall \gamma_{b}\}$,  
defines a group of diffeomorphisms of the fiber 
$F_{b}$, called the \textbf{holonomy} group of the assignment 
$ x \in M \mapsto (HT)_{x}\subset T_{x}M$, 
with reference to the point $b \in B$. (see Fig.~{\ref{fig4prop312}})
\end{prop312}

%(see Fig.~{\ref{fig4prop312}})

%    \include{mxrankbgarxiv.5}

\section{The Main Theorem}
\protect\label{ls5}

We begin this section with definitions of two new properties of maximal 
rank onto mappings: uniformly roughly isometric fibers 
[{\bf Definition~\ref{ldef55}}] 
and horizontal lift control [{\bf Definition~\ref{ldef56}}].

Let $\pi: M\rightarrow B$ be an onto smooth map with maximal rank 
between complete Riemannian manifolds $M$ and $B$ with dimensions 
$m$ and $n$, respectively. 

Consider $b_{0}\in B$ a fixed base point. 

For every $b_{1}, b_{2}\in B$ we will denote by 
$\gamma_{[b_{1}b_{2}]}$
a broken geodesic in $B$ joining $b_{1}$ to $b_{2}$. 
In particular $\gamma_{b_{1}}:= \gamma_{[b_{1}b_{1}]}$
will denote a broken geodesic loop at $b_{1}$.
Let $\varphi_{(\gamma_{[b_{1}b_{2}]})}: F_{b_{1}}
\longrightarrow F_{b_{2}}$ be the
corresponding O'Neill diffeomorphism to $\gamma_{[b_{1}b_{2}]}$, 
as in \textsf{Definition~\ref{ldef37}}.

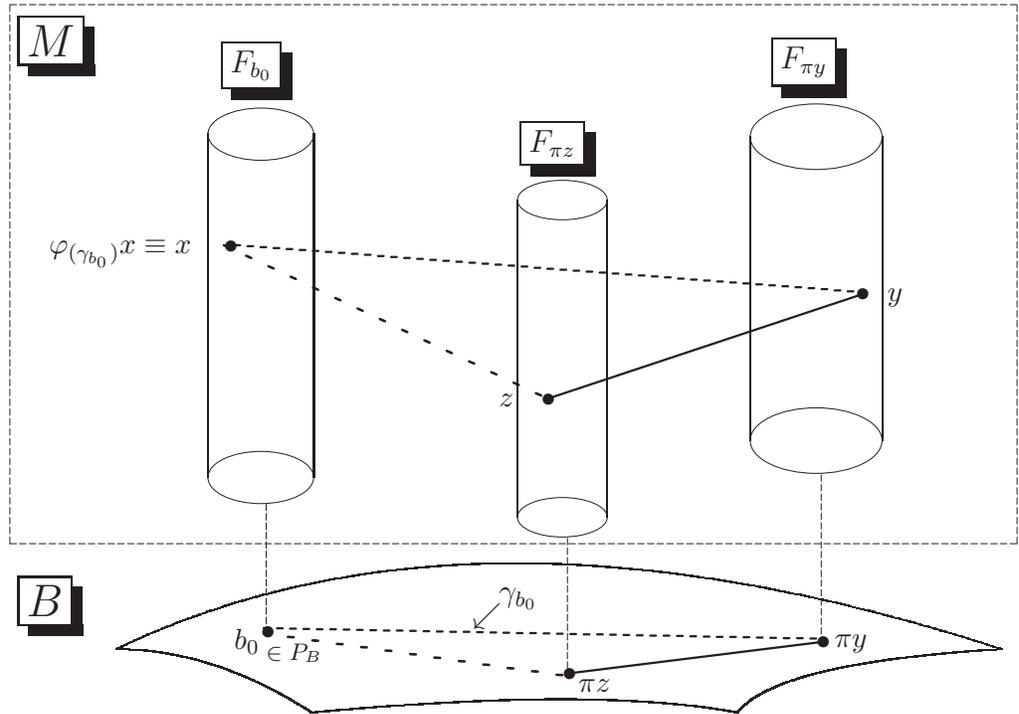
\begin{figure}[here]%oldthm

         \begin{picture}(390,270)(0,0)%(390,270)(0,0)

%\graphpaper[20](-75,0)(470,275)%default spacing is [10]
\definecolor{darkyellow}{rgb}{0.7,0.7,0}
\definecolor{orange}{rgb}{1.0,0.4,0}
\definecolor{darkgreen}{rgb}{0,0.8,0}
%\dottedline{2}(0,300)(390,300)%(0,270)(390,270)
%\dottedline{2}(0,0)(390,0)

\put(8,242){\shadowbox{\LARGE $M$}}
\dashline[+90]{3}(5,269)(386,269)
\dashline[+90]{3}(5,65)(5,269)
\dashline[+90]{3}(386,65)(386,269)
\dashline[+90]{3}(5,65)(386,65)

\put(85,233){\shadowbox{\large $F_{b_{0}}$}}

%cylinder
%parallels=horizontal curves
\put(100,220){\ellipse{40}{20}}
\put(100,90){\ellipse{40}{20}}
%meridians=vertical curves
\put(80,90){\line(0,1){130}}
\put(120,90){\line(0,1){130}}
\dashline[+90]{3}(102,80)(102,31)

\put(198,204){\shadowbox{\large $F_{\pi z}$}}

%cylinder
%parallels=horizontal curves
\put(214,195){\ellipse{34}{15}}
\put(214,75){\ellipse{34}{15}}
%meridians=vertical curves
\put(197,75){\line(0,1){120}}
\put(231,75){\line(0,1){120}}
\dashline[+90]{3}(216,67)(216,14)

\put(293,234){\shadowbox{\large $F_{\pi y}$}}

%cylinder
%parallels=horizontal curves
\put(310,219){\ellipse{50}{25}}
\put(310,104){\ellipse{50}{25}}
%meridians=vertical curves
\put(285,104){\line(0,1){115}}
\put(335,104){\line(0,1){115}}
\dashline[+90]{3}(312,91)(312,27)

\put(8,30){\shadowbox{\LARGE $B$}}
\qbezier(45,25)(180,90)(380,25)
\qbezier(119,1)(227,11)(280,1)
\qbezier(45,25)(82,25)(119,1)
\qbezier(380,25)(295,20)(280,1)

\thicklines
\put(217,16){\line(8,1){97}}%line pi z-pi y
%\put(217,16){{\color{orange}\line(8,1){97}}}%line pi z-pi y
\dashline{3}[0.7](101,31)(214,15)%dashline b-pi z
%%\put(217,16){\line(-7,1){115}}%line bo-pi z
%\put(217,16){{\color{green}\line(-7,1){115}}}%line bo-pi z
\dashline[+90]{3}(101,33)(310,29)%dashline b-pi y
%%\put(104,32){\line(45,-1){208}}%line bo-pi y
%\put(104,32){{\color{blue}\line(45,-1){208}}}%line bo-pi y
\thinlines

\put(90,25){$b_{0}$}
\put(102,21){\footnotesize $\in P_{B}$}
\put(100,29){$\bullet$}
\put(220,9){$\pi z$}
\put(214,13){$\bullet$}
\put(317,26){$\pi y$}
\put(310,25){$\bullet$}
\put(190,43){\large $\gamma_{b_{0}}$}
\put(180,35){$\swarrow$}

%text
\put(20,175){$\varphi_{({\gamma_{b_{0}}})}x\equiv x$}
\put(86,175){$\bullet$}
\thicklines
\dashline[+90]{3}(90,178)(329,160)
%{\color{blue}\dashline[+90]{3}(90,178)(329,160)}
\dashline{3}[0.7](86,178)(206,121)
%{\color{green}\dashline{3}[0.7](86,178)(206,121)}
\thinlines
\put(337,157){$y$}
\put(325,157){$\bullet$}
\put(190,117){$z$}
\put(206,117){$\bullet$}
\thicklines
\put(209,120){\line(3,1){120}}
%\put(209,120){{\color{orange}\line(3,1){120}}}
         \end{picture}

         \caption{trivial holonomy with reference to the point $b_{0}$.}
         \label{fig416.3thm53}
         \index{pictures!Theorem\ref{lthm53}}%oldthm
      \end{figure}

According to O"Neill~\cite{BON} an onto  maximal rank map  
$\pi: M\rightarrow B$ has \textbf{trivial holonomy} with reference 
to the point $b_{0}$, if for any 
broken geodesic loop $\gamma_{b_{0}}$, the corresponding O'Neill 
diffeomorphism $\varphi_{(\gamma_{b_{0}})}: F_{b_{0}} \longrightarrow
F_{b_{0}}$ is the identity map on $F_{b_{0}}$
(see Fig.~{\ref{fig416.3thm53}}).

%(see Fig.~{\ref{fig416.3thm53}}

\begin{def55} 
\protect\label{ldef55} 
An onto  maximal rank 
map  $\pi: M\rightarrow B$ has \textbf{uniformly} \textbf{roughly} 
\textbf{isometric} \textbf{fibers} 
\textbf{(RIF)}   
if for all $b \in B$  
there exist constants 
$A>1$ and $C>0$, both independent of $b$, such that,
\[
\frac{1}{A} d_{M}(x,x^{\prime}) - C \leq
d_{M}(\varphi_{(\gamma_{[b,b_{0}]})}(x),
\varphi_{(\gamma_{[b,b_{0}]})}(x^{\prime}))
\leq A \; d_{M}(x,x^{\prime}) + C
\]
for all $x,x^{\prime}\in F_{b}$, where $d_{M}$ denotes the Riemannian
metric on $M$.

In this case, since $\varphi_{(\gamma_{[b,b_{0}]})}$ is onto, it
follows that
$\varphi_{(\gamma_{[b,b_{0}]})}:F_{b}\longrightarrow F_{b_{0}}$ is a
rough isometry for each $b \in B$, and therefore the fibers are
uniformly roughly isometric (see Fig.~{\ref{fig416.4thm53}}).
\end{def55}

\begin{figure}[here]

         \begin{picture}(390,260)(0,0)%(390,270)(0,0)

%\graphpaper[20](-75,0)(470,275)%default spacing is [10]
\definecolor{darkyellow}{rgb}{0.7,0.7,0}
\definecolor{orange}{rgb}{1.0,0.4,0}
\definecolor{darkgreen}{rgb}{0,0.8,0}
%\dottedline{2}(0,300)(390,300)%(0,270)(390,270)
%\dottedline{2}(0,0)(390,0)

\put(8,242){\shadowbox{\LARGE $M$}}
\dashline[+90]{3}(5,269)(386,269)
\dashline[+90]{3}(5,65)(5,269)
\dashline[+90]{3}(386,65)(386,269)
\dashline[+90]{3}(5,65)(386,65)

\put(85,233){\shadowbox{\large $F_{b_{0}}$}}

%cylinder
%parallels=horizontal curves
\put(100,220){\ellipse{40}{20}}
\put(100,90){\ellipse{40}{20}}
%meridians=vertical curves
\put(80,90){\line(0,1){130}}
\put(120,90){\line(0,1){130}}
\dashline[+90]{3}(102,80)(102,31)

\put(293,234){\shadowbox{\large $F_{b}$}}

%cylinder
%parallels=horizontal curves
\put(310,219){\ellipse{50}{25}}
\put(310,104){\ellipse{50}{25}}
%meridians=vertical curves
\put(285,104){\line(0,1){115}}
\put(335,104){\line(0,1){115}}
\dashline[+90]{3}(312,91)(312,27)

\put(8,30){\shadowbox{\LARGE $B$}}
\qbezier(45,25)(180,90)(380,25)
\qbezier(119,1)(227,11)(280,1)
\qbezier(45,25)(82,25)(119,1)
\qbezier(380,25)(295,20)(280,1)

\thicklines
%\put(217,16){{\color{green}\line(-7,1){115}}}%line bo-pi q
\put(104,32){\line(45,-1){208}}%line bo-b
%\put(104,32){{\color{blue}\line(45,-1){208}}}%line bo-b
\thinlines

\put(97,20){$b_{0}\in P_{B}$}
\put(100,29){$\bullet$}
\put(317,26){$b\in P_{B}$}
\put(310,25){$\bullet$}
\put(200,36){$\gamma_{[b,b_{0}]}$}
%\put(200,36){$\color{blue}\gamma_{[b,b_{0}]}$}

%text
\put(30,175){$\varphi_{\left({\gamma_{[b,b_{0}]}}\right)}x$}
%\put(30,175){$\varphi_{\left({\color{blue}
%\gamma_{[b,b_{0}]}}\right)}x$}\put(86,175){$\bullet$}
\put(30,137)
  {$\varphi_{\left({\gamma_{[b,b_{0}]}}\right)}x^{\prime}$}
%\put(30,137)
%  {$\varphi_{\left({\color{blue}\gamma_{[b,b_{0}]}}\right)}x^{\prime}$}
\put(106,137){$\bullet$}
\dottedline[$\bullet$]{1}(108,141)(90,176)
\thicklines
{\dashline[+90]{3}(90,178)(329,160)}
%{\color{blue}\dashline[+90]{3}(90,178)(329,160)}
{\dashline{3}[0.7](108,140)(295,103)}
%{\color{blue}\dashline{3}[0.7](108,140)(295,103)}
\thinlines
\put(335,157){$x$}
\put(325,157){$\bullet$}
\put(335,100){$x^{\prime}$}
\put(295,100){$\bullet$}
\dottedline[$\bullet$]{1}(299,104)(328,159)

         \end{picture}

         \caption{\textbf{(RIF)}.}
         \label{fig416.4thm53}
         \index{pictures!Theorem\ref{lthm53}}
      \end{figure}

%(see Fig.~{\ref{fig416.4thm53}}

\begin{figure}[here]

         \begin{picture}(390,260)(0,0)%(390,270)(0,0)

%\graphpaper[20](-75,0)(470,275)%default spacing is [10]
\definecolor{darkyellow}{rgb}{0.7,0.7,0}
\definecolor{orange}{rgb}{1.0,0.4,0}
\definecolor{darkgreen}{rgb}{0,0.8,0}
%\dottedline{2}(0,300)(390,300)%(0,270)(390,270)
%\dottedline{2}(0,0)(390,0)

\put(8,242){\shadowbox{\LARGE $M$}}
\dashline[+90]{3}(5,269)(386,269)
\dashline[+90]{3}(5,65)(5,269)
\dashline[+90]{3}(386,65)(386,269)
\dashline[+90]{3}(5,65)(386,65)

\put(85,233){\shadowbox{\large $F_{b}$}}

%cylinder
%parallels=horizontal curves
\put(100,220){\ellipse{40}{20}}
\put(100,90){\ellipse{40}{20}}
%meridians=vertical curves
\put(80,90){\line(0,1){130}}
\put(120,90){\line(0,1){130}}
\dashline[+90]{3}(102,80)(102,31)

\put(293,234){\shadowbox{\large $F_{b^{\prime}}$}}

%cylinder
%parallels=horizontal curves
\put(310,219){\ellipse{50}{25}}
\put(310,104){\ellipse{50}{25}}
%meridians=vertical curves
\put(285,104){\line(0,1){115}}
\put(335,104){\line(0,1){115}}
\dashline[+90]{3}(312,91)(312,27)

\put(8,30){\shadowbox{\LARGE $B$}}
\qbezier(45,25)(180,90)(380,25)
\qbezier(119,1)(227,11)(280,1)
\qbezier(45,25)(82,25)(119,1)
\qbezier(380,25)(295,20)(280,1)

\put(93,29){$b$}
\put(100,29){$\bullet$}
\put(113,32){{$w$}}
%\put(113,32){{\color{blue}$w$}}
\put(317,26){$b^{\prime}$}
\put(310,25){$\bullet$}
\put(206,41){{$w^{\prime}$}}
%\put(206,41){{\color{green}$w^{\prime}$}}
\thicklines
\put(104,32){{\vector(4,-1){20}}}%vector w
%\put(104,32){{\color{blue}\vector(4,-1){20}}}%vector w
\put(310,28){{\vector(-4,1){112}}}%vector w'
%\put(310,28){{\color{green}\vector(-4,1){112}}}%vector w'
\thinlines

%text
\put(108,175){$x$}
\put(116,175){$\bullet$}
\put(145,181){{$v$}}
%\put(145,181){{\color{blue}$v$}}
\put(295,137){$x^{\prime}$}
\put(285,137){$\bullet$}
\put(230,142){{$v^{\prime}$}}
%\put(230,142){{\color{green}$v^{\prime}$}}
\thicklines
\put(120,178){{\vector(1,0){40}}}%vector v
%\put(120,178){{\color{blue}\vector(1,0){40}}}%vector v
{\dashline[+90]{3}(120,175)(226,175)}
%{\color{red}\dashline[+90]{3}(120,175)(226,175)}
\put(180,172){$\times$}
%\put(180,172){\color{red}$\times$}
\put(286,140){{\vector(-1,0){65}}}%vector v'
%\put(286,140){{\color{green}\vector(-1,0){65}}}%vector v'
{\dashline[+90]{3}(260,137)(286,137)}
%{\color{red}\dashline[+90]{3}(260,137)(286,137)}
\put(271,134){$\times$}
%\put(271,134){\color{red}$\times$}
\thinlines

         \end{picture}

         \caption{\textbf{(HLC)},  where $w$ %{\color{blue}$w$} 
         is short and $w^{\prime}$ %{\color{green}$w^{\prime}$} 
         is long.}
         \label{fig416.5thm53}
         \index{pictures!Theorem\ref{lthm53}}
      \end{figure}
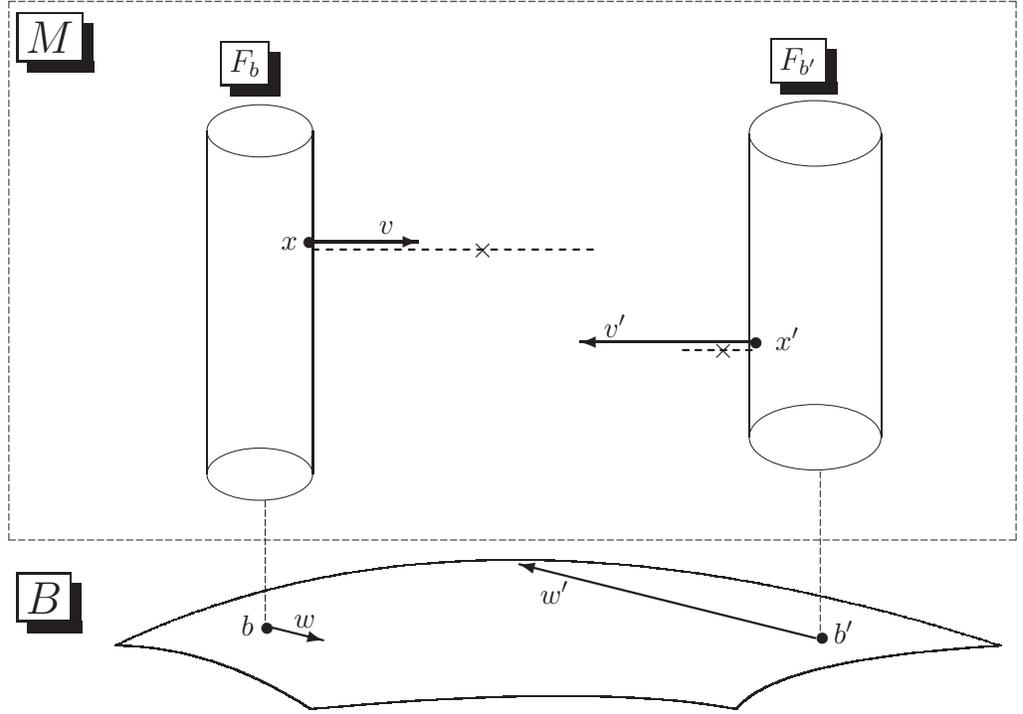

\begin{def56} 
\protect\label{ldef56} An onto  maximal rank map  $\pi: M\rightarrow B$ 
has \textbf{horizontal} \textbf{lift} \textbf{control} \textbf{(HLC)} if
for all $b \in B$ and for all
$x\in F_{b}$ there exist constants $\alpha\geq 1$ and $\beta >0$,
both independent of $b$ and $x$, such that
\[
\frac{1}{\alpha} ||w||_{B} - \beta \leq ||v||_{M} \leq \alpha
||w||_{B} + \beta
\]
for all $w\in T_{b}B$,  where $v$ is the unique horizontal lift of
$w$ through $x$ satisfying $||v||_{M}\leq 1$,
and $||\hspace{0.15in}||_{M}$, $||\hspace{0.15in}||_{B}$ denote the
inner product on $TM$ and $TB$, respectively 
(see Fig.~{\ref{fig416.5thm53}}).
\end{def56}

%(see Fig.~{\ref{fig416.5thm53}}

\vspace{0.1in}

Finally, we state and prove the main result. \textbf{Theorem~\ref{lthm53}} 
was motivated by O'Neill's~\cite{BON} question adapted for Mappings with 
Maximal Rank.

%\vspace{0.1in}

\begin{thm53}
\protect\label{lthm53}
Let $M$ and $B$ be complete Riemannian manifolds,  with 
bounded geometry and dimensions $m$ and $n$, respectively. 
Let $\pi: M\rightarrow B$ be an onto smooth maximal rank map, and 
let $b_{0}\in B$ be a fixed base point. 
Assume that $\pi$ has  
trivial holonomy,  uniformly roughly isometric fibers \textbf{(RIF)}, 
and horizontal lift control \textbf{(HLC)}. 
Then, $M$ is roughly isometric to the product $F_{b_{0}}\times B$.
\end{thm53}

In order to prove {\bf Theorem~\ref{lthm53}} we will need the
following technical Lemma and Proposition.

\begin{lem51}
\protect\label{llem51}
Let $A>1$, $C>0$, $\alpha\geq 1$, $\beta>0$, be given constants.
 
For any positive real numbers $\epsilon_{0}>0$, $\epsilon_{B}>0$ 
satisfying,
\begin{eqnarray}
\protect\label{leqn430}
\epsilon_{0} > C\cdot(1+A^{2}) & \mbox{ and } & \epsilon_{B} >
\left(\displaystyle{\frac{\epsilon_{0}-C}{A}}+\beta\right)\cdot\alpha
\end{eqnarray} 
the following hold:
\begin{equation}
\protect\label{lequ456}
      \left(\displaystyle{\frac{\epsilon_{0}-C}{A}}\right)>0
\end{equation}
\begin{equation}
\protect\label{lequ457}
      \displaystyle{\frac{\epsilon_{0}-C}{A}}<\epsilon_{0}
\end{equation}
\begin{equation}
\protect\label{lequ458}
      \displaystyle{\frac{\epsilon_{0}-C}{A^{2}}}-C>0
\end{equation}
\begin{equation}
\protect\label{lequ459}
      \displaystyle{\frac{\epsilon_{0}-C}{A}}>
       \displaystyle{\frac{\epsilon_{0}-C}{A^{2}}}-C
\end{equation}
\begin{equation}
\protect\label{lequ460}
      \displaystyle{\frac{1}{\alpha}\epsilon_{B}-\beta}>0
\end{equation}
\begin{equation}
\protect\label{lequ461}
      \displaystyle{\frac{\epsilon_{0}-C}{A}}<(\epsilon_{0}+C)\cdot A
\end{equation}
\begin{equation}
\protect\label{lequ462}
      2\left(\displaystyle{\frac{\epsilon_{0}-C}{A}}\right)<
      (2\epsilon_{0}+C)\cdot A
\end{equation}
\begin{equation}
\protect\label{lequ463}
      2\epsilon_{B}>2\left(
\displaystyle{\frac{\epsilon_{0}-C}{A}}+\beta\right)\alpha
\end{equation}
\end{lem51}
\pf\hspace{0.1in}

\underline{(\ref{lequ456})}: By the first inequality in (\ref{leqn430}),
\begin{eqnarray*}
\epsilon_{0} & > & C(\underbrace{1+A^{2}}_{>1})>C
\Rightarrow
\left(\displaystyle{\frac{\epsilon_{0}-C}{A}}\right)>0
\end{eqnarray*}

\underline{(\ref{lequ457})}: Since $A > 1$,
\begin{eqnarray*}
\epsilon_{0}(\underbrace{1-A}_{<0})<0<C & \Rightarrow &
\displaystyle{\frac{\epsilon_{0}}{A}}(1-A)
<\displaystyle{\frac{C}{A}}
\Rightarrow \\
& \Rightarrow &
\displaystyle{\frac{\epsilon_{0}}{A}}-\epsilon_{0}-
\displaystyle{\frac{C}{A}}<0
\Rightarrow%\\
%& \Rightarrow &
\displaystyle{\frac{\epsilon_{0}}{A}}-\displaystyle{\frac{C}{A}}<
\epsilon_{0} 
\Rightarrow \\
& \Rightarrow &
\displaystyle{\frac{\epsilon_{0}-C}{A}}<\epsilon_{0}
\end{eqnarray*}

\underline{(\ref{lequ458})}: By the first inequality in (\ref{leqn430}),
\begin{eqnarray*}
\epsilon_{0} > C(1+A^{2})
& \Rightarrow & %\epsilon_{0}>C+CA^{2} \Rightarrow
\epsilon_{0}-C-CA^{2}>0
\Rightarrow \displaystyle{\frac{\epsilon_{0}-C}{A^{2}}}-C>0
\end{eqnarray*}

\underline{(\ref{lequ459})}: Since $A > 1$,
\begin{eqnarray*}
& & 
\epsilon_{0}(\underbrace{1-A}_{<0})<0<C+CA 
(\underbrace{A-1}_{>0})
\Rightarrow\\ 
& \Rightarrow & \epsilon_{0}(1-A)<C+CA^{2}-CA \Rightarrow %\\
\epsilon_{0}-C-CA^{2}<\epsilon_{0}A-CA
\Rightarrow \\
& & \\
& \Rightarrow & \displaystyle{\frac{\epsilon_{0}-C}{A^{2}}}-C<
\displaystyle{\frac{\epsilon_{0}-C}{A}}
\end{eqnarray*}

\underline{(\ref{lequ460})}: By the second inequality in (\ref{leqn430}),
\begin{eqnarray*}
& & \epsilon_{B}  > 
\left(\displaystyle{\frac{\epsilon_{0}-C}{A}}+\beta\right)\alpha=
\underbrace{\left(\displaystyle{\frac{\epsilon_{0}-C}{A}}
          \right)\alpha}_{>0} + \beta\alpha > \beta\alpha\Rightarrow\\ 
& \stackrel{\div \alpha}{\Rightarrow} & 
\displaystyle{\frac{1}{\alpha}\epsilon_{B}-\beta}>0
\end{eqnarray*}

\underline{(\ref{lequ461})}: $A>1, C>0$ and $\epsilon_{0}>0$ imply that,
\begin{eqnarray*}
& & \epsilon_{0}(1-A)(1+A)<0 \Rightarrow\\
& & \\
& \Rightarrow & \epsilon_{0}(1-A)(1+A)<0<
C(1+A^{2})\Rightarrow 
\epsilon_{0}(1-A^{2})<C(1+A^{2}) \Rightarrow\\
& & \\
& \Rightarrow &
\epsilon_{0}-\epsilon_{0}A^{2}<C+CA^{2} \Rightarrow 
\epsilon_{0}-C<\epsilon_{0}A^{2}+CA^{2} \Rightarrow \\
& & \\
& \Rightarrow &
\displaystyle{\frac{\epsilon_{0}-C}{A}}<(\epsilon_{0}+C)A
\end{eqnarray*}

\underline{(\ref{lequ462})}: From $A>1, C>0$ and $\epsilon_{0}>0$ 
we have,
\begin{eqnarray*}
\lefteqn{\left\{\begin{array}{l}
2\epsilon_{0}(1-A)(1+A)<0\\
0<C(2+A^{2})
\end{array}\right.%}
\Rightarrow 2\epsilon_{0}(1-A)(1+A)<0<2C+CA^{2}\Rightarrow}\\
& & \\
& \Rightarrow & 2\epsilon_{0}(1-A^{2})<2C+CA^{2}\Rightarrow
2\epsilon_{0}-2\epsilon_{0}A^{2}<2C+CA^{2} \Rightarrow \\
& & \\
& \Rightarrow & 2\epsilon_{0}-2C<2\epsilon_{0}A^{2}+CA^{2}\Rightarrow \\
& & \\
& \Rightarrow &
2\left(\displaystyle{\frac{\epsilon_{0}-C}{A}}\right)<(2\epsilon_{0}+C)A
\end{eqnarray*}

\underline{(\ref{lequ463})}: By the second inequality in (\ref{leqn430}),
\begin{eqnarray*}
& & 
\epsilon_{B} > \left(\displaystyle{\frac{\epsilon_{0}-C}{A}}\right)
\alpha+\beta\alpha>\left(\displaystyle{\frac{\epsilon_{0}-C}{A}}\right)
\alpha+\displaystyle{\frac{\beta\alpha}{2}}\Rightarrow\\
& \stackrel{\times 2}{\Rightarrow} &
2\epsilon_{B} > 2\left(\displaystyle{\frac{\epsilon_{0}-C}{A}}\right)
\alpha+\beta\alpha  \Rightarrow \\%=
& \Rightarrow & 2\epsilon_{B} >  
\left[2\left(\displaystyle{\frac{\epsilon_{0}-C}{A}}\right)
+\beta\right]\alpha
\end{eqnarray*}

\pfe

\begin{prop52}
\protect\label{lprop52}
Suppose that $M$ and $B$  are complete 
$m-$dimensional and $n-$dimensional 
Riemannian ma\-ni\-folds, respectively, 
both with bounded geometry. %(see \textsf{Definition~\ref{ldef24}}).
Let $\pi : M\rightarrow B$ be an onto smooth map with maximal rank,
and $b_{0}\in B$ be fixed. 

Assume that  $\{\phi_{b} : F_{b} \rightarrow F_{b_{0}}\}_{b\in B}$
is a family of bijective rough isometries satisfying,
\begin{eqnarray}
\protect\label{leqn464}
\lefteqn{\forall b \in B, \exists   A>1, \exists C>0 :} \nonumber\\
& & \frac{1}{A} d_{M}(x,x^{\prime}) - C \leq
d_{M}(\phi_{b}(x),\phi_{b}(x^{\prime}))
\leq A d_{M}(x,x^{\prime}) + C, \\
& & \forall x,x^{\prime}\in F_{b} \nonumber 
\end{eqnarray}
where, $A$ and $C$ are universal constants independent of $b$.

Then, the following hold:
\begin{itemize}
\item If $P_{0}$ is an $\epsilon_{0}-$separated set and
$\epsilon_{0}-$full in $F_{b_{0}}$, where 
we assume that $\epsilon_{0}>C$, then the set
\[
P_{b}:= \phi_{b}^{-1}(P_{0})
\]
is an $\hat{\epsilon}-$separated set and $\tilde{\epsilon}-$full in
$F_{b}$, where
$\hat{\epsilon}:= \left(\displaystyle{\frac{\epsilon_{0}-C}{A}}
\right)>0$ and $\tilde{\epsilon}:= (\epsilon_{0}+C)A>0$.
\item For all $b\in B$ the corresponding nets $P_{b}$ are
\underline{uniformly} roughly isometric to $P_{0}$ with respect to the
combinatorial metric $\delta$.
\end{itemize}
\end{prop52}
\pf\hspace{0.1in}
 From Lemma~2.5~\cite{MK1}, 
 %in~\cite{MK1}, on page 397 %M. Kanai reference(16),
a complete Riemannian manifold with bounded geometry is roughly
isometric to each of its nets. 

This implies that for each $b\in B$,
\[
(P_{b}, \delta) \stackrel{{\bf R.I.}}{\longrightarrow} (F_{b}, d_{M})
\Rightarrow
\]
\begin{equation}
\protect\label{lequ427}
\frac{1}{2\hat{\epsilon}} d_{M}(p_{1},p_{2}) \leq
\delta (p_{1},p_{2}) \leq \tilde{a} d_{M}(p_{1},p_{2}) + \tilde{c},
\hspace{0.2in} \forall p_{1},p_{2}\in P_{b}%\hspace{0.2in} {\bf [1]}
\end{equation}
where $\tilde{a}:=\tilde{a}(m, k_{M}, \tilde{\epsilon})>1,
       \tilde{c}:=\tilde{c}(m, k_{M}, \tilde{\epsilon})>0$, and
       $(F_{b}, d_{M})$ indicates that on each fiber $F_{b}$ we will
       use the induced Riemannian metric from $M$.

\vspace{0.1in}

Also, by~\cite{MK1} (Lemma~5) we have,
\[
(P_{0}, \delta_{0}) \stackrel{{\bf R.I.}}{\longrightarrow}
(F_{0}, d_{M})
\Rightarrow
\]
\[
\frac{1}{2{\epsilon_{0}}} d_{M}(p_{3},p_{4}) \leq
\delta_{0} (p_{3},p_{4}) \leq \tilde{a}_{0} d_{M}(p_{3},p_{4}) +
\tilde{c}_{0} \Rightarrow
\]
\begin{equation}
\protect\label{lequ428}
\frac{1}{\tilde{a}_{0}} \delta_{0} (p_{3},p_{4}) -
\frac{\tilde{c}_{0}}{\tilde{a}_{0}}\leq
d_{M}(p_{3},p_{4}) \leq 2{\epsilon_{0}} \delta_{0} (p_{3},p_{4}),
\hspace{0.2in}\forall p_{3},p_{4}\in P_{0} %\hspace{0.2in} {\bf [2]}
\end{equation}
where $\tilde{a}_{0}:=\tilde{a}_{0}(m, k_{M}, \epsilon_{0})>1,
       \tilde{c}_{0}:=\tilde{c}_{0}(m, k_{M}, \epsilon_{0})>0$, and
       $(F_{0}, d_{M})$ indicates that on the fiber $F_{0}$
       the induced Riemannian metric from $M$ is used.

\vspace{0.1in}

From  (\ref{leqm464}), for all $p_{1},p_{2}\in F_{b}$,
\[
(F_{b}, d_{M}) \stackrel{{\bf R.I.}}{\longrightarrow} (F_{0}, d_{M})
\Rightarrow
\]
\begin{equation}
\protect\label{lequ429}
\frac{1}{A} d_{M}(\phi_{b}(p_{1}),\phi_{b}(p_{2})) - \frac{C}{A} \leq
d_{M}(p_{1},p_{2})
\leq A d_{M}(\phi_{b}(p_{1}),\phi_{b}(p_{2})) + AC
%\hspace{0.2in} {\bf [3]}
\end{equation}

Next, we observe the following diagram for $\iota=1,2$,

\[
\begin{array}{ccccccc}
P_{b} & & F_{b} & & F_{0} & & P_{0} \\
p_{\iota} & \hookrightarrow  & p_{\iota} &
\stackrel{\phi_{b}}{\rightarrow} & \phi_{b}(p_{\iota}) &
\hookrightarrow  & \phi_{b}(p_{\iota})
\end{array}
\]
where, $p_{\iota} \in P_{b}:= \phi_{b}^{-1}{P_{0}}\Rightarrow
\phi_{b}(p_{\iota}) \in P_{0}$.

\vspace{0.1in}

We claim that,
\[
(P_{b}, \delta) \stackrel{{\bf unif. R.I.}}{\longrightarrow}
(P_{0}, \delta_{0})
\]

\vspace{0.1in}

Indeed, let $p_{1},p_{2}\in P_{b}$.

\vspace{0.1in}

By  (\ref{lequ427}), (\ref{lequ429}) and (\ref{lequ428}), we may write,
\begin{eqnarray*}
& \stackrel{(\ref{lequ427})}{\Longrightarrow} &
\frac{1}{2\hat{\epsilon}} d_{M}(p_{1},p_{2}) \leq
\delta (p_{1},p_{2}) \leq \tilde{a} d_{M}(p_{1},p_{2}) + \tilde{c}
\Rightarrow\\
& & \\
& \stackrel{(\ref{lequ429})}{\Longrightarrow} &
\frac{1}{2\hat{\epsilon} A} d_{M}(\phi_{b}(p_{1}),\phi_{b}(p_{2}))
- \frac{C}{2\hat{\epsilon }A}\leq
\delta (p_{1},p_{2}) =\\
& & = \delta (p_{1},p_{2}) \leq
\tilde{a} A
d_{M}(\phi_{b}(p_{1}),\phi_{b}(p_{2})) + \tilde{a} A C + \tilde{c}
\Rightarrow\\
& & \\
& \stackrel{(\ref{lequ428})}{\Longrightarrow} &
\frac{1}{2\hat{\epsilon} A\tilde{a}_{0}}
\delta_{0} (\phi_{b}(p_{1}),\phi_{b}(p_{2})) -
\frac{\tilde{c}_{0}}{2\hat{\epsilon} A \tilde{a}_{0}} -
\frac{C}{2\hat{\epsilon} A}\leq \delta (p_{1},p_{2}) =\\
& & \\
& & = \delta (p_{1},p_{2}) \leq
\tilde{a}_{0} A 2{\epsilon_{0}}
\delta_{0} (\phi_{b}(p_{1}),\phi_{b}(p_{2}))
\tilde{a} A C + \tilde{c}%\\
\end{eqnarray*}
which can be rewritten as,
\[
\frac{1}{A_{net}}
\delta_{0} (\phi_{b}(p_{1}),\phi_{b}(p_{2})) - C_{net}
\leq \delta (p_{1},p_{2})\leq
A_{net} \delta_{0} (\phi_{b}(p_{1}),\phi_{b}(p_{2})) + C_{net}
\]
where,
\begin{eqnarray*}
A_{net} & :=  & A_{net}(m, k_{M}, \epsilon_{0}, C, A):=
2A\max\{\hat{\epsilon}\tilde{a}_{0}, \tilde{a}\epsilon_{0}\}\geq 1\\
C_{net} & := & C_{net}(m, k_{M}, \epsilon_{0}, C, A):=
\max\left\{\tilde{a} A C + \tilde{c},
\frac{1}{2\hat{\epsilon} A\tilde{a}_{0}}\left(
\frac{\tilde{c}_{0}}{\tilde{a}_{0}} + C\right)\right\}>0
\end{eqnarray*}
and the Proposition is proved.

\pfe

%PROOF OF THEOREM

\textit{Proof of Theorem~\ref{lthm53}}. 

In order to prove the Theorem, by~\cite{MK1} (Lemma~2.5),
%{\bf Lemma~\ref{llem24}}, 
it suffices to show that an $\epsilon-$net in $M$ is roughly 
isometric to an $\epsilon'-$net in the product $F_{b_{0}}\times B$. 
We remark here that in the proof of~\cite{MK1} (Lemma~5) 
the maximal property of an $\epsilon-$net is not required, it 
sufficient that the "net" be a countable, $\epsilon-$separated 
and $\epsilon-$full set.

\vspace{0.1in}

We will proceed with the proof by constructing in 2 steps a rough
isometry $\phi$ between countable, separated full sets in $M$ and in
$F_{b_{0}}\times B$.

\vspace{0.1in}

In {\it Step 1.} we combine the diffeomorphisms
$\varphi_{(\gamma_{[b,b_{0}]})}$ with two countable maximal separated
sets, $P_{0}$ in the fiber $F_{b_{0}}\subset M$ and $P_{B}$ in $B$, in a
fashion that will produce a suitable countable separated full set $P$
in $M$. We also show that the product $P_{0}\times P_{B}$ is a
countable separated full set in $F_{b_{0}}\times B$.

Then, in {\it Step 2.} we introduce a bijection $\phi$ from $P$ to
$P_{0}\times P_{B}$, which will turn out to be the rough isometry
between discrete approximations of $M$ and $F_{b_{0}}\times B$, as
mentioned above.

\underline{\it Step 1.}

Let the positive constants $A,C$ and $\alpha, \beta$, be as in 
conditions ({\bf RIF}) and ({\bf HLC}),  respectively.
Let us choose and fix two constants $\epsilon_{0}>0$ and 
$\epsilon_{B}>0$ satisfying the inequalities (see (\ref{leqn430})),
\[
\epsilon_{0} > C\cdot(1+A^{2})\hspace{0.5in}\mbox{ and }\hspace{0.5in}
\epsilon_{B} >
\left(\displaystyle{\frac{\epsilon_{0}-C}{A}}+\beta\right)\cdot\alpha
\]

We first define two countable sets $P_{0}\subseteq F_{b_{0}}\subset M$
and $P_{B}\subseteq B$, with $b_{0}\in P_{B}$, where $P_{0}$ is a
maximal $\epsilon_{0}-$separated set,
\[
\forall p,q \in P_{0}, p\neq q \Rightarrow d_{M}(p,q)\geq \epsilon_{0}
\]
and $P_{B}$ is a maximal $\epsilon_{B}-$separated set,
\[
\forall b_{1},b_{2} \in P_{B}, b_{1}\neq b_{2} \Rightarrow
d_{B}(b_{1},b_{2})\geq \epsilon_{B}
\]
and then we introduce  the net structure
$N_{0}=\{N_{0}(p): p\in P_{0}\}$  of $P_{0}$ given by,
\[
N_{0}(p)=\{q\in P_{0}: 0<d_{M}(p,q)\leq 2\epsilon_{0}\}
\]
and  $N_{B}=\{N_{B}(b): b\in P_{B}\}$ the net structure of
$P_{B}$ defined by,
\[
N_{B}(b)=\{\hat{b}\in P_{B}: 0<d_{B}(b,\hat{b})\leq 2\epsilon_{B}\}
\]

Observe that {\bf Proposition~\ref{lprop23}} implies $P_{0}$ is
$\epsilon_{0}-$full in $F_{b_{0}}$ and $P_{B}$ is $\epsilon_{B}-$full
in $B$.

\vspace{0.1in}

We now, construct $P$ a countable
$\left(\displaystyle{\frac{\epsilon_{0}-C}{A}}\right)-$separated full
set in $M$.

\vspace{0.1in}

For each $b\in P_{B}$, let us look first at
$\varphi_{(\gamma_{[b,b_{0}]})}^{-1}(P_{0})\subseteq F_{b}$.

We claim that,
\begin{equation}
\protect\label{lequ432}
\varphi_{(\gamma_{[b,b_{0}]})}^{-1}(P_{0})\hspace{0.1in}
\mbox{ is a countable }
\left(\displaystyle{\frac{\epsilon_{0}-C}{A}}\right)-\mbox{separated
subset of }F_{b}
\end{equation}

The set $\varphi_{(\gamma_{[b,b_{0}]})}^{-1}(P_{0})$ is countable, due
to the fact that $P_{0}$ is countable and
$\varphi_{(\gamma_{[b,b_{0}]})}$ is bijective.

Notice that by~(\ref{lequ456}), 
$\left(\displaystyle{\frac{\epsilon_{0}-C}{A}}\right)>0$
%(see {\bf Lemma~\ref{llem51} (1.)})

\vspace{0.1in}

Now, for $b\in P_{B}$ let us consider either $b=b_{0}$ or $b\neq b_{0}$.

\vspace{0.1in}

If \underline{$b=b_{0}$}, since $P_{0}$ is $\epsilon_{0}-$separated, 
by (\ref{lequ457}) we have  
$\displaystyle{\frac{\epsilon_{0}-C}{A}}<\epsilon_{0}$
%(see {\bf Lemma~\ref{llem51} (2.)})
, so we conclude that
$\varphi_{(\gamma_{[b,b_{0}]})}^{-1}\left(P_{0}\right)$ is
$\left(\displaystyle{\frac{\epsilon_{0}-C}{A}}\right)-$separated, which
is claim (\ref{lequ432}).

If \underline{$b\neq b_{0}$}, because $\varphi_{(\gamma_{[b,b_{0}]})}$ 
is a diffeomorphism, we have,
\begin{eqnarray*}
\lefteqn{\forall p,q\in \varphi_{(\gamma_{[b,b_{0}]})}^{-1}(P_{0}):
p\neq q \Rightarrow}\\
& \Rightarrow & \varphi_{(\gamma_{[b,b_{0}]})}(p)\neq
\varphi_{(\gamma_{[b,b_{0}]})}(q)
\mbox{ in }P_{0} \Rightarrow
\epsilon_{0}\leq d_{M}(\varphi_{(\gamma_{[b,b_{0}]})}(p),
\varphi_{(\gamma_{[b,b_{0}]})}(q)) \leq\\
& \stackrel{(\bf RIF)}{\leq} & A d_{M}(p,q)+C \Rightarrow
d_{M}(p,q)\geq
\left(\displaystyle{\frac{\epsilon_{0}-C}{A}}\right)
\end{eqnarray*}

and claim (\ref{lequ432}) follows.

\vspace{0.1in}

\begin{figure}[here]

         \begin{picture}(370,270)(0,0)

%\dottedline{2}(0,286)(370,286)
%\dottedline{2}(0,0)(370,0)

\put(330,253){\shadowbox{\Huge $M$}}
\dashline[+90]{3}(20,283)(368,283)
\dashline[+90]{3}(20,85)(20,283)
\dashline[+90]{3}(368,85)(368,283)
\dashline[+90]{3}(20,85)(368,85)

\put(45,245){\shadowbox{\large $F_{b_{0}}$}}
%cylinder
%parallels=horizontal curves
\put(59,240){\ellipse{18}{6}}
\put(59,110){\ellipse{18}{6}}

%meridians=vertical curves
\put(50,110){\line(0,1){130}}
\put(68,110){\line(0,1){130}}

%\thicklines
%\put(100,110){\line(-1,3){38}}
%\thinlines
\put(51,112){$\bullet$}
\put(35,112){$p_{l}$}
%\dashline[+90]{3}(100,110)(307,110)
\put(57,125){$\bullet$}
%\dashline[+90]{3}(62,123)(250,123)
\put(35,135){$\vdots$}
\put(62,145){$\bullet$}
%\dashline[+90]{3}(94,143)(295,143)
\put(56,160){$\bullet$}
\put(35,160){$p_{3}$}
%\dashline[+90]{3}(67,158)(285,158)
\put(65,185){$\bullet$}
\put(35,185){$p_{2}$}
%\dashline[+90]{3}(104,183)(265,183)
\put(53,205){$\bullet$}
\put(35,205){$p_{1}$}
%\dashline[+90]{3}(86,203)(275,203)
\put(64,226){$\bullet$}
\put(35,226){$p_{0}$}
%\dashline[+90]{3}(62,224)(269,224)
\put(30,95){$\in P_{0}$}

%text
\put(170,268){$\varphi^{-1}_{(\gamma_{[b,b_{0}]})}$}
\thicklines
\put(75,260){\vector(1,0){215}}
\thinlines

\put(294,245){\shadowbox{\large $F_{b}$}}
%cylinder
%parallels=horizontal curves
\put(307,240){\ellipse{12}{4}}
\put(307,110){\ellipse{12}{4}}

%meridians=vertical curves
\put(301,110){\line(0,1){130}}
\put(313,110){\line(0,1){130}}

%\thicklines
%\put(196,110){\line(2,3){77}}
%\thinlines
\put(305,112){$\bullet$}
\put(317,112){$\varphi^{-1}_{(\gamma_{[b,b_{0}]})}p_{l}$}
\put(308,125){$\bullet$}
\put(335,135){$\vdots$}
\put(302,145){$\bullet$}
\put(309,160){$\bullet$}
\put(317,160){$\varphi^{-1}_{(\gamma_{[b,b_{0}]})}p_{3}$}
\put(311,185){$\bullet$}
\put(317,185){$\varphi^{-1}_{(\gamma_{[b,b_{0}]})}p_{2}$}
\put(301,205){$\bullet$}
\put(317,205){$\varphi^{-1}_{(\gamma_{[b,b_{0}]})}p_{1}$}
\put(299,226){$\bullet$}
\put(317,226){$\varphi^{-1}_{(\gamma_{[b,b_{0}]})}p_{0}$}
\put(324,95){$\in P$}

%cylinder
%parallels=horizontal curves
\put(153,240){\ellipse{6}{2}}
\put(153,110){\ellipse{6}{2}}

%meridians=vertical curves
\put(156,110){\line(0,1){130}}
\put(150,110){\line(0,1){130}}

\put(147,112){$\bullet$}
\put(152,125){$\bullet$}
\put(148,145){$\bullet$}
\put(153,160){$\bullet$}
\put(151,185){$\bullet$}
\put(149,205){$\bullet$}
\put(150,226){$\bullet$}
\put(155,95){$\in P$}

%cylinder
%parallels=horizontal curves
\put(218,240){\ellipse{24}{8}}
\put(218,110){\ellipse{24}{8}}

%meridians=vertical curves
\put(206,110){\line(0,1){130}}
\put(230,110){\line(0,1){130}}

\put(211,112){$\bullet$}
\put(227,125){$\bullet$}
\put(218,145){$\bullet$}
\put(207,160){$\bullet$}
\put(215,185){$\bullet$}
\put(203,205){$\bullet$}
\put(219,226){$\bullet$}
\put(220,95){$\in P$}

\put(330,34){\shadowbox{\Huge $B$}}
\qbezier(33,25)(186,90)(340,25)
\qbezier(83,5)(186,20)(290,5)
\qbezier(33,25)(65,25)(83,5)
\qbezier(290,5)(300,25)(340,25)

\thicklines
\qbezier(59,25)(150,5)(304,29)
\thinlines
\put(165,27){\ovalbox{\large $\gamma_{[b,b_{0}]}$}}
\put(201,18){\large $\leadsto$}
\put(56,23){$\bullet$}
\put(62,27){$b_{0}$}
\dashline[+90]{3}(59,29)(59,108)
\put(304,27){$\bullet$}
\put(310,24){$b$}
\dashline[+90]{3}(307,29)(307,108)
\put(150,37){$\bullet$}
\dashline[+90]{3}(153,39)(153,108)
\put(215,47){$\bullet$}
\dashline[+90]{3}(218,49)(218,108)
\put(310,10){$\in P_{B}$}

%\put(198,100){$q\in P$}
%\put(105,60){\shadowbox{\scriptsize minimal geodesics}}
%\put(128,78){\line(0,1){52}}
%\put(128,130){\vector(-1,0){33}}

%\put(157,78){\line(0,1){52}}
%\put(157,130){\vector(1,0){50}}

%\put(179,70){\line(1,0){60}}
%\put(239,70){\vector(0,-1){38}}
   \end{picture}

         \caption{The net
   $P=\bigcup_{b\in P_{B}}\varphi_{(\gamma_{[b,b_{0}]})}^{-1}(P_{0})$.}
         \label{fig415thm53}
         \index{pictures!Theorem\ref{lthm53}}
      \end{figure}
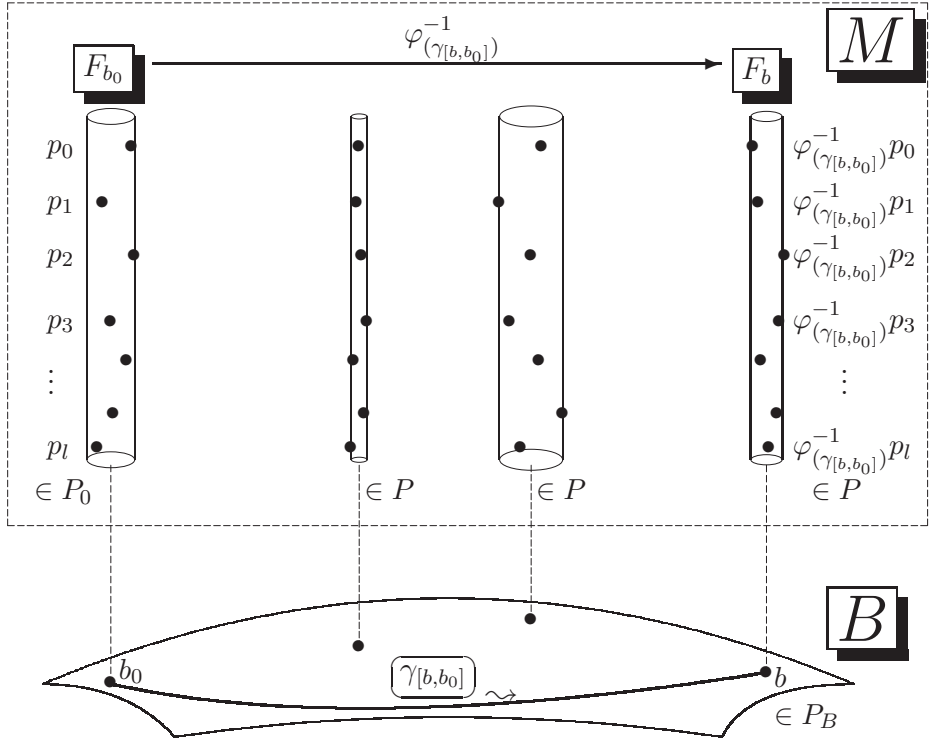

Let (see Fig.~{\ref{fig415thm53}}),
\[
P:=\bigcup_{b\in P_{B}}\varphi_{(\gamma_{[b,b_{0}]})}^{-1}(P_{0})
\]

%(see Fig.~{\ref{fig415thm53}})

Since $P_{0}, P_{B}$ are countable sets and
$\varphi_{(\gamma_{[b,b_{0}]})}$
is a bijection for all $b\in B$, the set $P$ is also
\underline{countable}.

\vspace{0.1in}

To show that $P$ is
$\left(\displaystyle{\frac{\epsilon_{0}-C}{A}}\right)-$separated, we
proceed as follows.

\vspace{0.1in}

For any $p,q\in P$ such that $p\neq q$, we have only two cases,

{\it (CASE I:)} $\exists b\in P_{B}: p,q\in
\varphi_{(\gamma_{[b,b_{0}]})}^{-1}(P_{0}).$

In that case (\ref{lequ432}) gives us,
$d_{M}(p,q)\geq \displaystyle{\frac{\epsilon_{0}-C}{A}}.$

{\it (CASE II:)} $\exists b\in P_{B}: p\in
\varphi_{(\gamma_{[b,b_{0}]})}^{-1}(P_{0})$ and
$\exists \tilde{b}\in P_{B}: \linebreak
q\in \varphi_{(\gamma_{[\tilde{b},b_{0}]})}^{-1}(P_{0})$, where
$\pi p=b\neq \tilde{b}=\pi\tilde{b}.$

We claim,
\[
d_{M}(p,q)\geq \displaystyle{\frac{1}{\alpha}\epsilon_{B}-\beta}
\stackrel{(\ref{leqn430})}{>} \displaystyle{\frac{\epsilon_{0}-C}{A}}
\]
 We will only verify the first inequality, since the second one is
the requirement (\ref{leqn430}) on $\epsilon_{B}$.

Let $\varsigma$ be a general curve parametrized proportionally to .a.l. curve joining $p$ and $q$ in
$M$, with length $\ell(\varsigma)$. In this case $\pi\circ\varsigma$ is
a curve joining $b$ and $\tilde{b}$. In addition, we will denote by
$\gamma_{b\tilde{b}}$ the minimal geodesic joining $b$ and $\tilde{b}$.

For all $x\in M$, we can write
$v=v_{V}\oplus v_{H}\in T_{x}M=(VT)_{x}\oplus (HT)_{x}$ in a 
unique way, where %by {\bf Lemma~\ref{llem43}} 
$v_{V}\in (\ker(\pi_{\ast})_{x})$ and
$v_{H}\in (\ker(\pi_{\ast})_{x})^{\perp}$. Furthermore,
$||v||_{M}=||v_{V}\oplus v_{H}||_{M}=
\sqrt{||v_{V}||_{M}^{2}+||v_{H}||_{M}^{2}}.$

By property {(\bf HLC)}, the facts that $\gamma_{b\tilde{b}}$
is a minimal geodesic and $B$ is complete,
\begin{eqnarray*}
\ell(\varsigma)&=& \int_{0}^{1} ||\varsigma'(t)||_{M}  dt =
\int_{0}^{1} ||\varsigma'_{H}(t)\oplus \varsigma'_{V}(t)||_{M} dt
\geq \int_{0}^{1} ||\varsigma'_{H}(t)||_{M} dt
\stackrel{(\bf HLC)}{\geq} \\
&\geq & \displaystyle{\frac{1}{\alpha}}
\int_{0}^{1}||(\pi\circ\varsigma)'(t)||_{B} dt -\beta=
\displaystyle{\frac{1}{\alpha}} \ell(\pi\circ\varsigma)-\beta
\stackrel{min.geod.}{\geq}\\
&\geq & \displaystyle{\frac{1}{\alpha}}\ell(\gamma_{b\tilde{b}})- \beta
\stackrel{B complete}{=}
\displaystyle{\frac{1}{\alpha}}d_{B}(b,\tilde{b})- \beta,
\hspace{0.2in}\forall \varsigma
\end{eqnarray*}
which is a lower bound on the length of any curve $\varsigma$ joining
$p$ and $q$ in $M$, independent of the curve $\varsigma$.

Finally, by definition of infimum and from $b,\tilde{b}\in P_{B}$,
\begin{eqnarray*}
d_{M}(p,q):=\inf_{\varsigma\subset M}\ell(\varsigma)\geq
\displaystyle{\frac{1}{\alpha}}d_{B}(b,\tilde{b})- \beta
\stackrel{b,\tilde{b}\in P_{B}}{\geq}
\displaystyle{\frac{1}{\alpha}\epsilon_{B}-\beta}
\stackrel{(\ref{lequ460})%Lem~\ref{llem51}(5.)
}{>}0
\end{eqnarray*}
and the claim is proved.

So, $P$ is
$\left(\displaystyle{\frac{\epsilon_{0}-C}{A}}\right)-$
\underline{separated}.

We introduce the net structure $N_{P}=\{N_{P}(p): p\in P\}$  of $P$
given by,
\[
N_{P}(p)=\left\{q\in P: 0<d_{M}(p,q)\leq
2\left(\displaystyle{\frac{\epsilon_{0}-C}{A}}\right)\right\}
\]

Next, we prove that that $P$ is
$[(\epsilon_{0}+C)A+\alpha\epsilon_{B}+\beta]-$full in $M$, i.e.
\[
M=B_{[(\epsilon_{0}+C)A+\alpha\epsilon_{B}+\beta]}P=
\{x\in M: d_{M}(x,P)<(\epsilon_{0}+C)A+\alpha\epsilon_{B}+\beta\}
\]

We want to show that,
\[
d_{M}(x,P):=\inf_{p\in P}d_{M}(x,p)<
(\epsilon_{0}+C)A+\alpha\epsilon_{B}+\beta,
\hspace{0.2in}\forall x\in M
\]

Let $x\in M$. Either $x\in P$ or $x\in M\setminus P$.

If $x\in P$, then $d_{M}(x,P)=0<(\epsilon_{0}+C)A+\alpha\epsilon_{B}
+\beta$.

If $x\in M\setminus P$, let $b:=\pi x\in B$.

There are two  cases for such $b$,
either $b\in P_{B}$ or $b\in B\setminus P_{B}$.

{\it (CASE I:)} If $b\in P_{B}$, since $x$ is not in $P$,
$P_{0}$ is maximal, property {(\bf RIF)} holds, and $b\in P_{B}$
implies
$\varphi_{(\gamma_{[b,b_{0}]})}^{-1}(p_{0})\in P,
\forall p_{0}\in P_{0}$, then,

\begin{eqnarray*}
b\in P_{B} & \stackrel{x\in M\setminus P}{\Longrightarrow}&
\varphi_{(\gamma_{[b,b_{0}]})}x \in F_{b_{0}}\setminus P_{0}
\stackrel{P_{0} max.}{\Longrightarrow} \exists p_{0}\in P_{0}:
d_{M}(\varphi_{(\gamma_{[b,b_{0}]})}x,p_{0})<\epsilon_{0}\\
& \stackrel{(\bf RIF)}{\Longrightarrow} &
\displaystyle{\frac{1}{A}}
d_{M}(x,\varphi_{(\gamma_{[b,b_{0}]})}^{-1}(p_{0}))-C
\leq d_{M}(\varphi_{(\gamma_{[b,b_{0}]})}x,p_{0})<
\epsilon_{0}\Rightarrow \\
&\Rightarrow &
d_{M}(x,\underbrace{\varphi_{(\gamma_{[b,b_{0}]})}^{-1}(p_{0})}_{\in P})
< (\epsilon_{0}+C)A<(\epsilon_{0}+C)A+
\underbrace{\alpha\epsilon_{B}+\beta}_{>0}\Rightarrow \\
&\Rightarrow &
\inf_{p\in P}d_{M}(x,p)\leq
d_{M}(x,\varphi_{(\gamma_{[b,b_{0}]})}^{-1}(p_{0}))
< (\epsilon_{0}+C)A+\alpha\epsilon_{B}+\beta\\
&\Rightarrow & d_{M}(x,P)<(\epsilon_{0}+C)A+\alpha\epsilon_{B}+\beta
\end{eqnarray*}

{\it (CASE II:)} If $b\in B\setminus P_{B}$, by the maximality
of $P_{B}$ there exists
$\bar{b}\in P_{B}: d_{B}(b,\bar{b})<\epsilon_{B}$.

We wish to obtain $\bar{p}\in P$ satisfying
$d_{M}(x,\bar{p})\leq (\epsilon_{0}+C)A+\alpha\epsilon_{B}+\beta$,
which will be accomplished as follows.

Let,
\begin{eqnarray*}
\gamma_{b\bar{b}}: [t_{1},t_{2}] & \longrightarrow & B,  
\gamma_{b\bar{b}}(t_{1}):=b, \gamma_{b\bar{b}}(t_{2}):=\bar{b}
\end{eqnarray*}
be a parametrization proportional to arclength of a  
minimal geodesic joining $b$ and
$\bar{b}$ in $B$, and let $\Gamma_{b\bar{b}}$ be its unique horizontal
lift through $x$, which in particular satisfies
\begin{eqnarray}
\protect\label{leqn455} 
\Gamma_{b\bar{b}}(t_{2}) & \in  & F_{\bar{b}}
\end{eqnarray}

We have, by (\ref{leqn455}), the fact that $P_{0}$ is 
$\epsilon_{0}-$full
in $F_{b_{0}}$ and the definition of infimum,
\begin{eqnarray}
\protect\label{leqn433}
& & \Gamma_{b\bar{b}}(t_{2})\in F_{\bar{b}} \Rightarrow
\varphi_{(\gamma_{[\bar{b},b_{0}]})}\left(\Gamma_{b\bar{b}}(t_{2})
\right) \in F_{b_{0}}\Longrightarrow \nonumber \\
%\stackrel{Prop\ref{lprop23}}{\Longrightarrow}
& \stackrel{P_{0}-full}{\Longrightarrow} &
d_{M}\left(\varphi_{(\gamma_{[\bar{b},b_{0}]})}
\left(\Gamma_{b\bar{b}}(t_{2})
\right),P_{0}\right)<\epsilon_{0}\Rightarrow \nonumber\\
& \stackrel{infimum}{\Longrightarrow} &
\exists p_{0}\in P_{0}:d_{M}\left(\varphi_{(\gamma_{[\bar{b},b_{0}]})}
\left(\Gamma_{b\bar{b}}(t_{2})\right), p_{0}\right)<\epsilon_{0}
\end{eqnarray}

We claim that the desired $\bar{p}\in P$ is exactly
$\varphi_{(\gamma_{[\bar{b},b_{0}]})}^{-1}(p_{0})$.

Indeed, 
\[
\bar{b}\in P_{B}, \hspace{0.1in} 
\varphi_{(\gamma_{[\bar{b},b_{0}]})}^{-1}
(p_{0})\in \varphi_{(\gamma_{[\bar{b},b_{0}]})}^{-1}(P_{0}) \Rightarrow
\varphi_{(\gamma_{[\bar{b},b_{0}]})}^{-1}(p_{0})\in P
\]

Furthermore, by the triangle inequality,
\begin{equation}
\protect\label{lequ434}
d_{M}\left(\varphi_{(\gamma_{[\bar{b},b_{0}]})}^{-1}(p_{0}),x\right)\leq
d_{M}\left(\varphi_{(\gamma_{[\bar{b},b_{0}]})}^{-1}(p_{0}),
\Gamma_{b\bar{b}}(t_{2})\right)
+ d_{M}\left(\Gamma_{b\bar{b}}(t_{2}),x\right)
\end{equation}

In addition we have,
\begin{eqnarray}
\protect\label{leqn435}
& & \displaystyle{\frac{1}{A}}d_{M}\left(\Gamma_{b\bar{b}}(t_{2}),
\varphi_{(\gamma_{[\bar{b},b_{0}]})}^{-1}(p_{0})\right)-C
\stackrel{(\bf RIF)}{\leq}d_{M}
\left(\varphi_{(\gamma_{[\bar{b},b_{0}]})}
\left(\Gamma_{b\bar{b}}(t_{2})\right), p_{0}\right) \Rightarrow
\nonumber\\
&\Rightarrow &
d_{M}\left(\Gamma_{b\bar{b}}(t_{2}),
\varphi_{(\gamma_{[\bar{b},b_{0}]})}^{-1}(p_{0})\right) \leq
\left[d_{M}\left(\varphi_{(\gamma_{[\bar{b},b_{0}]})}
\left(\Gamma_{b\bar{b}}(t_{2})\right), p_{0}\right)+C\right]A <
\nonumber\\
&  & \stackrel{(\ref{leqn433})}{<}(\epsilon_{0}+C)A
\end{eqnarray}

Also,
\begin{eqnarray}
\protect\label{leqn436}
\lefteqn{d_{M}\left(\Gamma_{b\bar{b}}(t_{2}),x\right) =
d_{M}\left(\Gamma_{b\bar{b}}(t_{2}),\Gamma_{b\bar{b}}(t_{1})\right)
\stackrel{infimum}{\leq} \ell(\Gamma_{b\bar{b}})=
\int_{0}^{1}||\Gamma'_{b\bar{b}}(t)||_{M}  dt}
\nonumber\\
& \stackrel{(\bf HLC)}{\leq} & \alpha\int_{0}^{1}
||\left(\pi\circ\Gamma_{b\bar{b}}\right)'(t)||_{B} dt + \beta
\stackrel{{\bf (hl.1)}}{=}\alpha\int_{0}^{1}
||\gamma'_{b\bar{b}}(t)||_{B} dt + \beta = \nonumber\\
& = & \alpha \ell(\gamma_{b\bar{b}})
+ \beta \stackrel{min.geod.}{=} \alpha d_{B}(b,\bar{b}) + \beta
\stackrel{def.\bar{b}}{<} \alpha\epsilon_{B} + \beta
\end{eqnarray}

Thus, by combining (\ref{lequ434}), (\ref{leqn435}) and
(\ref{leqn436}),
\[
d_{M}\left(\varphi_{(\gamma_{[\bar{b},b_{0}]})}^{-1}(p_{0}),x\right)
\leq
(\epsilon_{0}+C)A + \alpha\epsilon_{B} + \beta
\]
which in turn implies,
\[
d_{M}(x,P)=\inf_{p\in P}d_{M}(x,p)\leq
d_{M}(x,\underbrace{
\varphi_{(\gamma_{[\bar{b},b_{0}]})}^{-1}(p_{0})}_{\in P})
\leq (\epsilon_{0}+C)A + \alpha\epsilon_{B} + \beta
\]
and we conclude that $P$ is
$[(\epsilon_{0}+C)A+\alpha\epsilon_{B}+\beta]-$
\underline{full in $M$}.

In what follows, we show that $P_{0}\times P_{B}$ is countable,
$(\epsilon_{0}+\epsilon_{B})$-separated, and
%by {\bf Proposition~\ref{lprop23}} $P_{0}\times P_{B}$ is
$(\epsilon_{0}+\epsilon_{B})$-full in $F_{b_{0}}\times B$.

In $F_{b_{0}}\times B$ we have the induced product metric from
$M\times B$,
\[
d_{\times}\left((x,b),(\tilde{x},\tilde{b})\right):=
d_{M}(x,\tilde{x}) + d_{B}(b,\tilde{b})
\]
for all $x,\tilde{x}\in F_{b_{0}}$ and $b,\tilde{b}\in B$.

$P_{0}\times P_{B}$ being \underline{countable} comes from the fact
that both $P_{0}$ and $P_{B}$ have that property.

Let $(x,b)\neq (\tilde{x},\tilde{b})\in P_{0}\times P_{B}$. Since
$P_{0}$ is $\epsilon_{0}-$separated and $P_{B}$ is
$\epsilon_{B}-$separated,
\[
d_{\times}\left((x,b),(\tilde{x},\tilde{b})\right)=
d_{M}(x,\tilde{x}) + d_{B}(b,\tilde{b})\geq \epsilon_{0} + \epsilon_{B}
\]
and $P_{0}\times P_{B}$ is
$(\epsilon_{0}+\epsilon_{B})$-\underline{separated}.

To prove that $P_{0}\times P_{B}$ is
$(\epsilon_{0}+\epsilon_{B})$-full in $F_{b_{0}}\times B$, i.e.,
\[
F_{b_{0}}\times B=\{(x,b)\in F_{b_{0}}\times B:
d_{\times}\left((x,b), P_{0}\times P_{B}\right)<\epsilon_{0} +
\epsilon_{B}\}
\] we need to show that,
\[
d_{\times}\left((x,b), P_{0}\times P_{B}\right)<\epsilon_{0} +
\epsilon_{B}, \hspace{0.2in}\forall (x,b)\in F_{b_{0}}\times B
\]

Let $(x,b)\in F_{b_{0}}\times B$.

Since, $P_{0}$ is $\epsilon_{0}$-full in $F_{b_{0}}$, there exists
$p_{0}\in P_{0}: d_{M}(p_{0},x)<\epsilon_{0}$. Similarly, $P_{B}$
being $\epsilon_{B}$-full in $B$, implies that there exists
$\bar{b}\in B: d_{B}(\bar{b},b)<\epsilon_{B}$.

Therefore,
\begin{eqnarray*}
d_{\times}\left((x,b), P_{0}\times P_{B}\right)&
\stackrel{inf.}{\leq} & d_{\times}\left((x,b), (p_{0},\bar{b})\right)
=d_{M}(x,p_{0}) + d_{B}(b,\bar{b})<\\
& < & \epsilon_{0} + \epsilon_{B}
\end{eqnarray*}
and since $(x,b)$ is arbitrary, we conclude that $P_{0}\times P_{B}$
is $(\epsilon_{0}+\epsilon_{B})$-\underline{full} in
$F_{b_{0}}\times B$.

\underline{\it Step 2.}

Let us initially define some notation as well as provide a geometric
interpretation of a "net". We will assume that all nets are connected,
since we can repeat the argument on each connected component of the
underlying manifold.

\begin{figure}[here]

         \begin{picture}(310,130)(0,0)

%\dottedline{2}(0,140)(310,140)
%\dottedline{2}(0,0)(310,0)

\put(50,80){\circle{40}}
\put(48,77){$\bullet$}
\put(45,86){$p_{0}$}
\dashline[+90]{3}(50,80)(63,23)
\put(63,23){\circle{40}}
\put(61,20){$\bullet$}
\put(58,14){$p_{1}$}
\dashline[+90]{3}(63,23)(120,115)
\put(120,115){\circle{40}}
\put(118,112){$\bullet$}
\put(118,122){$\updownarrow$}
\put(124,122){$\frac{\epsilon}{2}$}
\dashline[+90]{3}(120,115)(155,90)
\put(155,90){\circle{40}}
\put(153,87){$\bullet$}
\put(148,81){$p_{\imath-1}$}
\dashline[+90]{3}(155,90)(210,110)
\put(185,50){\vector(0,1){50}}
\put(136,36){\doublebox{$%\epsilon\leq
             d(p_{\imath-1},p_{\imath})%\leq 2\epsilon
             \in [\epsilon,2\epsilon]$}}
\put(210,110){\circle{40}}
\put(208,107){$\bullet$}
\put(205,101){$p_{\imath}$}
\dashline[+90]{3}(210,110)(282,107)
\put(282,107){\circle{40}}
\put(280,104){$\bullet$}
\put(278,98){$p_{\ell}$}
         \end{picture}

      \caption{A discrete path $(p_{0},p_{1},\ldots,p_{\ell})$ in an
      $\epsilon-$net, and its elements regarded as centers of coins
      with diameter $\epsilon$.}
         \label{fig416thm53}
         \index{pictures!Theorem\ref{lthm53}}
      \end{figure}
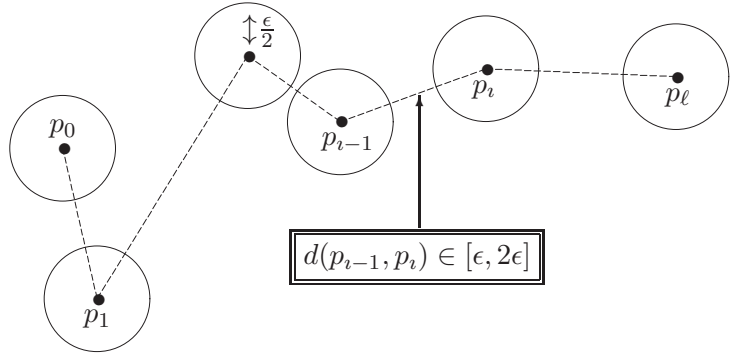

Graphically, we will describe an element of an $\epsilon-$net as the
center of a ball of radius $\frac{\epsilon}{2}$, which can be
visualized as a coin with diameter $\epsilon$. So, the control of
distances between elements in an $\epsilon-$net allows us to
describe it as a countable set of coins, which can touch boundaries
but can never overlap. We will call such element an
$\mathbf{\epsilon-}$\textbf{coin} (see Fig.~{\ref{fig416thm53}}).

%(see Fig.~{\ref{fig416thm53}})

We define a map $\phi$ from $P\subset M$ into
$P_{0}\times P_{B}\subset F_{b_{0}}\times B$ as follows,
\[
\begin{array}{rrcl}
\phi: & P=\bigcup_{b\in P_{B}}
\varphi_{(\gamma_{[b,b_{0}]})}^{-1}(P_{0})&
\longrightarrow & P_{0}\times P_{B}\\
& p & \longmapsto &
\phi p:=\left(\varphi_{(\gamma_{[\pi p,b_{0}]})}p,\pi p\right)
\end{array}
\]

\vspace{0.1in}

\cl{1} $\phi$ is well-defined.

\vspace{0.1in}

If $p\in \varphi_{(\gamma_{b_{0}})}^{-1}(P_{0})=P_{0}$, then
$\pi p=b_{0}$ and
$\phi p=\left(\varphi_{(\gamma_{b_{0}})}p,b_{0}\right)=
(p,b_{0})\in P_{0}\times P_{B}$.

If $p\in \varphi_{(\gamma_{[b,b_{0}]})}^{-1}(P_{0})$, where $b\in P_{B},
b\neq b_{0}$, then $\pi p=b$ and
$\phi p=\left(\varphi_{(\gamma_{[b,b_{0}]})}p,b\right)\in P_{0}
\times P_{B}$.

\vspace{0.1in}

\cl{2} $\phi$ is 1-1.

\vspace{0.1in}

Let $p,q\in P$ and assume that $\phi p=\phi q$.

Thus,
\[
\phi p=\phi q\Rightarrow
\left\{
\begin{array}{l}
\varphi_{(\gamma_{[\pi p,b_{0}]})}p = \varphi_{(\gamma_{[\pi q,b_{0}]})}q\\
\pi p=\pi q=:b
\end{array}\right.%}
\Rightarrow
\varphi_{(\gamma_{[b,b_{0}]})}p = \varphi_{(\gamma_{[b,b_{0}]})}q
\stackrel{diffeo.}\Rightarrow p=q
\]and $\phi$ is injective.

\vspace{0.1in}

\cl{3} $\phi$ is onto.

\vspace{0.1in}

Let $(p_{0},b)\in P_{0}\times P_{B}$ and define $p\in F_{b}$ by,
\[
p:=
\varphi_{(\gamma_{[b,b_{0}]})}^{-1}p_{0}\in
\varphi_{(\gamma_{[b,b_{0}]})}^{-1}(P_{0})
\subseteq F_{b}
%\subseteq F_{b} %\hspace{0.1in}
\] where, $\pi p=b$ and $p=p_{0}$ if $b=b_{0}$.

Thus,
\[
\phi p=\left(\varphi_{(\gamma_{[\pi p,b_{0}]})}p,\pi p\right)=
\left(\varphi_{(\gamma_{[b,b_{0}]})}p,b\right)=(p_{0},b)
\] and $\phi$ is surjective.

\vspace{0.1in}

\cl{4} $\phi$ satisfies {\bf (RI.1)}.

\vspace{0.1in}

By (Claim 3), for any given $\epsilon >0$, we have,
\[
P_{0}\times P_{B}=\phi\left(P\right)=
{\cal B}_{\epsilon}\phi\left(P\right)=\{(p_{0},b)\in P_{0}\times P_{B}:
d_{\times}\left((p_{0},b), \phi\left(P\right)\right)<\epsilon\}
\] in other words, $\phi$ is $\epsilon-$full in $P_{0}\times P_{B}$
for any $\epsilon >0$, which is exactly {\bf (RI.1)} for $\phi$.

\vspace{0.1in}

\cl{5} $\phi$ satisfies {\bf (RI.2)}.

\vspace{0.1in}

We want to show that there exist constants $a\geq 1$ and $c>0$
satisfying,
\[
\displaystyle{\frac{1}{a}}\delta_{P}(p,q)-c\leq
\delta_{\times}(\phi p, \phi q)\leq a\delta_{P}(p,q) + c,
\hspace{0.2in} \forall p,q\in P
\] where $\delta_{P}$ is the combinatorial metric of $P$, and
$\delta_{\times}$ is the discrete product metric of $P_{0}\times P_{B}$
given by,
\[
\delta_{\times}((p,b),(\tilde{p},\tilde{b})):=
\delta_{P_{0}}(p,\tilde{p}) + \delta_{P_{B}}(b,\tilde{b})
\]
for all $(p,b),(\tilde{p},\tilde{b}) \in P_{0}\times P_{B}$.

In terms of $\delta_{P_{0}}$ and $\delta_{P_{B}}$, the condition we
want to verify for $\phi$, translates into,
$\exists a\geq 1, \exists c>0: \forall p,q\in P$,
\begin{equation}
\protect\label{lequ437}
%\lefteqn{\exists a\geq 1, c>0: \forall p,q\in P}\\
\displaystyle{\frac{1}{a}}\delta_{P}(p,q)-c\leq
\delta_{P_{0}}
\left(\varphi_{(\gamma_{[\pi p,b_{0}]})}p,
      \varphi_{(\gamma_{[\pi q,b_{0}]})}q\right)
+ \delta_{P_{B}}(\pi p,\pi q) \leq a\delta_{P}(p,q) + c
\end{equation}
%for all $p,q\in P$.

\vspace{0.1in}

Indeed, let $p,q \in P=\bigcup_{b\in P_{B}}
\varphi_{(\gamma_{[b,b_{0}]})}^{-1}(P_{0})$.

%\vspace{0.1in}

Let $\gamma_{\mbox{\scriptsize min}}$ be a minimal geodesic joining
$\pi_{q}$ to $\pi_{p}$ in $B$, and let its unique horizontal lift
through $q$ be parametrized by
$\Gamma_{q}: [t_{1}, t_{2}]\longrightarrow M$.

    \begin{figure}[here]%oldthm

         \begin{picture}(390,270)(0,0)

%\graphpaper[20](-75,0)(470,275)%default spacing is [10]
\definecolor{darkyellow}{rgb}{0.7,0.7,0}
\definecolor{orange}{rgb}{1.0,0.4,0}
\definecolor{darkgreen}{rgb}{0,0.8,0}
%\dottedline{2}(0,270)(390,270)
%\dottedline{2}(0,0)(390,0)

\put(8,242){\shadowbox{\LARGE $M$}}
\dashline[+90]{3}(5,269)(386,269)
\dashline[+90]{3}(5,65)(5,269)
\dashline[+90]{3}(386,65)(386,269)
\dashline[+90]{3}(5,65)(386,65)

\put(85,233){\shadowbox{\large $F_{b_{0}}$}}

%cylinder
%parallels=horizontal curves
\put(100,220){\ellipse{40}{20}}
\put(100,90){\ellipse{40}{20}}
%meridians=vertical curves
\put(80,90){\line(0,1){130}}
\put(120,90){\line(0,1){130}}
\dashline[+90]{3}(102,80)(102,31)

\put(198,204){\shadowbox{\large $F_{\pi q}$}}

%cylinder
%parallels=horizontal curves
\put(214,195){\ellipse{34}{15}}
\put(214,75){\ellipse{34}{15}}
%meridians=vertical curves
\put(197,75){\line(0,1){120}}
\put(231,75){\line(0,1){120}}
\dashline[+90]{3}(216,67)(216,14)

\put(293,234){\shadowbox{\large $F_{\pi p}$}}

%cylinder
%parallels=horizontal curves
\put(310,219){\ellipse{50}{25}}
\put(310,104){\ellipse{50}{25}}
%meridians=vertical curves
\put(285,104){\line(0,1){115}}
\put(335,104){\line(0,1){115}}
\dashline[+90]{3}(312,91)(312,27)

\put(8,30){\shadowbox{\LARGE $B$}}
\qbezier(45,25)(180,90)(380,25)
\qbezier(119,1)(227,11)(280,1)
\qbezier(45,25)(82,25)(119,1)
\qbezier(380,25)(295,20)(280,1)

\thicklines
\put(217,16){\line(8,1){97}}%line pi q-pi p
\dashline{3}[0.7](101,31)(215,15)%dashline bo-pi q
%%\put(217,16){\line(-7,1){115}}%line bo-pi q
\dashline[+90]{3}(101,32)(311,28)%line bo-pi p
%%\put(104,32){\line(45,-1){208}}%line bo-pi p
\thinlines

\put(90,25){$b_{0}$}
\put(102,21){\footnotesize $\in P_{B}$}
\put(100,29){$\bullet$}
\put(201,10){$\pi q$}
\put(217,8){\footnotesize $\in P_{B}$}
\put(214,13){$\bullet$}
%\dashline[+90]{3}(196,29)(196,109)%connects q-pi q
\put(317,26){$\pi p$}
\put(331,26){\footnotesize $\in P_{B}$}
\put(310,25){$\bullet$}
%\dashline[+90]{3}(271,29)(271,223)%connects p-pi p
\put(178,38){{\large $\mathbf\gamma_{[\pi p,b_{0}]}$}}
%\put(160,35){$\swarrow$}
\put(139,16){{\large $\mathbf\gamma_{[\pi q,b_{0}]}$}}
%\put(135,13){$\nearrow$}
\put(257,11){{\large
              $\mathbf\gamma_{\mbox{\scriptsize min}}$}}
%\put(248,13){$\nwarrow$}

%text
\put(200,248){{\large $\mathbf \varphi_{(\gamma_{[\pi p,b_{0}]})}$}}
\thicklines
\put(290,240){\vector(-1,0){173}}
\thinlines
\put(139,214){{\large $\mathbf \varphi_{(\gamma_{[\pi q,b_{0}]})}$}}
\thicklines
\put(196,216){\vector(-4,1){79}}
\put(241,218){{\large $\mathbf 
               \varphi_{(\gamma_{\mbox{\scriptsize min}})}$}}
\thicklines
\put(228,222){\vector(4,1){63}}
\thinlines

\put(30,186){$\varphi_{(\gamma_{[\pi p,b_{0}]})}p$}
\put(50,170){$\in P_{0}$}
\put(110,186){$\bullet$}
\thicklines
\dashline[+90]{3}(112,189)(291,197)
\thinlines
\put(14,119){\footnotesize $\varphi_{(\gamma_{[\pi p,b_{0}]})}
             \Gamma_{q}(t_{2})$}
\put(86,110){$\bullet$}
\put(18,105){$\equiv\varphi_{(\gamma_{[\pi q,b_{0}]})}q$}
\put(50,88){$\in P_{0}$}
\thicklines
\dashline{3}[0.7](88,113)(206,76)%\dashline[+90]{3}(68,113)(186,76)
\dashline[+90]{3}(88,113)(327,105)
\thinlines

\put(338,194){$p\in P$}
\put(291,194){$\bullet$}
\put(336,110){$\varphi_{(\gamma_{\mbox{\scriptsize min}})}q$}
\put(336,94){$\equiv\Gamma_{q}(t_{2})$}
\put(325,102){$\bullet$}
\put(166,70){$q\in P$}
\put(206,72){$\bullet$}
\dottedline[$\diamond$]{3}(212,78)(291,193)%connects p-q

\put(255,96){{\large $\mathbf \Gamma_{q}$}}
\thicklines
\put(208,75){\line(4,1){120}}%connects q-\Gamma_{q}(t_{2}) 

%\put(308,103){\line(-1,4){24}}
\dottedline[$\bullet$]{1}(294,196)(328,105)%connects p-\Gamma_{q}(t_{2})
\dottedline[$\bullet$]{1}(112,188)(89,112)%connects images of 
                                   %p and \Gamma_{q}(t_{2}) in F_{b_{0}}
         \end{picture}

         \caption{An upper bound for the distance
         $d_{M}(p,q)$, via the triangle inequality.}
         \label{fig416.1thm53}
         \index{pictures!Theorem\ref{lthm53}}
      \end{figure}
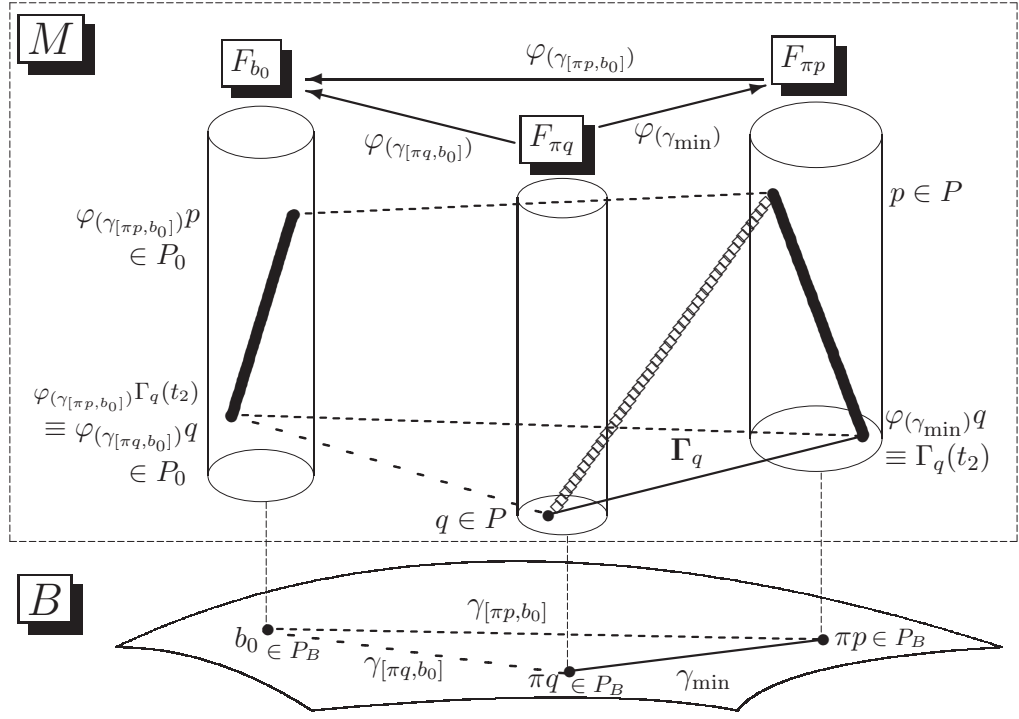

%(see Fig.~{\ref{fig416.1thm53}})

By the triangle inequality, the definition of distance,
{\bf Proposition~\ref{lprop35}}, trivial holonomy {\bf (TH)}, 
and property {\bf (RIF)}, we have (see Fig.~{\ref{fig416.1thm53}})
\begin{eqnarray*}
\lefteqn{d_{M}(p,q)  \stackrel{\triangle}{\leq}
d_{M}(p,\varphi_{(\gamma_{\mbox{\scriptsize min}})}q%\Gamma_{q}(t_{2})
) +  d_{M}(\varphi_{(\gamma_{\mbox{\scriptsize min}})}q%\Gamma_{q}(t_{2})
,q)
\leq} \nonumber\\
&  \stackrel{dist.}{\leq} & 
d_{M}(p,\varphi_{(\gamma_{\mbox{\scriptsize min}})}q%\Gamma_{q}(t_{2})
) + \ell(\Gamma_{q})
\leq \nonumber\\
&  \stackrel{(RIF)}{\leq} &  
A d_{M}(\varphi_{(\gamma_{[\pi p,b_{0}]})}p,
\varphi_{(\gamma_{[\pi p,b_{0}]})}\varphi_{(\gamma_{\mbox{\scriptsize min}})}q%\Gamma_{q}(t_{2})
) + AC +
\ell(\Gamma_{q})
\leq \nonumber\\
& \stackrel{Prop~\ref{lprop35}}{\leq} &
A d_{M}(\varphi_{(\gamma_{[\pi p,b_{0}]})}p,
          \varphi_{(\gamma_{[\pi p,b_{0}]})}\varphi_{(\gamma_{\mbox{\scriptsize min}})}q%\Gamma_{q}(t_{2})
) + AC
+ \alpha
\left[\ell(\gamma_{\mbox{\scriptsize min}})+\beta (t_{2}-t_{1})\right]
\nonumber\\
 & \stackrel{\bf (TH)}{\leq} &
A d_{M}(\varphi_{(\gamma_{[\pi p,b_{0}]})}p,
          \varphi_{(\gamma_{[\pi q,b_{0}]})}q) + AC
+ \alpha
\left[\ell(\gamma_{\mbox{\scriptsize min}})+\beta (t_{2}-t_{1})\right]
\nonumber\\
& \stackrel{dist.}{=} &
\alpha d_{B}(\pi p, \pi q) +
A d_{M}(\varphi_{(\gamma_{[\pi p,b_{0}]})}p,
\varphi_{(\gamma_{[\pi q,b_{0}]})}q) + \alpha\beta (t_{2}-t_{1}) + AC 
\nonumber\\
\end{eqnarray*}
i.e.,
\begin{equation}
\protect\label{lequ438}
d_{M}(p,q) \leq \alpha d_{B}(\pi p, \pi q) +
A d_{M}(\varphi_{(\gamma_{[\pi p,b_{0}]})}p,
\varphi_{(\gamma_{[\pi q,b_{0}]})}q) + AC + 
\alpha\beta (t_{2}-t_{1})
\end{equation}

Now, by Lemma~2.5 in~\cite{MK1} we have for the nets $P$, $P_{0}$ and 
$P_{B}$, respectively,
\begin{equation}
\protect\label{lequ439}
\exists \hat{a}(m,k_{M},\hat{\epsilon}) \geq 1,
\exists \hat{c}(m,k_{M},\hat{\epsilon}) > 0: \;
\frac{1}{\hat{a}}\delta_{P}(p,q) - \frac{\hat{c}}{\hat{a}}
\leq d_{M}(p,q)
\end{equation}
\begin{equation}
\protect\label{lequ440}
d_{M}(\varphi_{(\gamma_{[\pi p,b_{0}]})}p,
\varphi_{(\gamma_{[\pi q,b_{0}]})}q)\leq
2\epsilon_{0}
\delta_{P_{0}}(\varphi_{(\gamma_{[\pi p,b_{0}]})}p,
\varphi_{(\gamma_{[\pi q,b_{0}]})}q)
\end{equation}
\begin{equation}
\protect\label{lequ441}
d_{B}(\pi p, \pi q) \leq 2\epsilon_{B}\delta_{P_{B}}(\pi p, \pi q)
\end{equation}
where $m=\dim M$, $k_{M}>0$ [see \ref{ldef24}(BG.R)] and 
$\hat{\epsilon}>0$ [see \ref{lprop35}] are constants. 

If we combine (\ref{lequ439}), (\ref{lequ440}),  and
(\ref{lequ441}) into (\ref{lequ438}), we obtain,
\begin{eqnarray}
\protect\label{leqn442}
& & \frac{1}{\hat{a}}\delta_{P}(p,q) - \frac{\hat{c}}{\hat{a}}
\stackrel{(\ref{lequ439})}{\leq} d_{M}(p,q)
\leq \nonumber\\
& & \stackrel{(\ref{lequ438})}{\leq} 
\alpha d_{B}(\pi p, \pi q) +
A d_{M}(\varphi_{(\gamma_{[\pi p,b_{0}]})}p,
\varphi_{(\gamma_{[\pi q,b_{0}]})}q) +  AC +
\alpha\beta (t_{2}-t_{1})
\nonumber\\
& & \stackrel{(\ref{lequ440})}{\leq} 
\alpha d_{B}(\pi p, \pi q) + A\left(2\epsilon_{0}
\delta_{P_{0}}(\varphi_{(\gamma_{[\pi p,b_{0}]})}p,
\varphi_{(\gamma_{[\pi q,b_{0}]})}q)
\right) +  AC + 
\nonumber\\
& & \hspace{0.4in} + \alpha\beta (t_{2}-t_{1}) \leq 
\nonumber\\
& & \stackrel{(\ref{lequ441})}{\leq} 
\alpha 2\epsilon_{B}\delta_{P_{B}}(\pi p, \pi q) +
A\left(2\epsilon_{0}
\delta_{P_{0}}(\varphi_{(\gamma_{[\pi p,b_{0}]})}p,
\varphi_{(\gamma_{[\pi q,b_{0}]})}q)
\right) +  AC + \nonumber \\
& & \hspace{0.4in} + \alpha\beta (t_{2}-t_{1}) \Rightarrow 
\nonumber\\
& \Rightarrow &
\frac{1}{\hat{a}}\delta_{P}(p,q) - \frac{\hat{c}}{\hat{a}}
\leq 2A\epsilon_{0}
\delta_{P_{0}}(\varphi_{(\gamma_{[\pi p,b_{0}]})}p,
\varphi_{(\gamma_{[\pi q,b_{0}]})}q)
+ 2\alpha\epsilon_{B}\delta_{P_{B}}(\pi p, \pi q) + \nonumber \\
& &  +  AC + \alpha\beta (t_{2}-t_{1}) \Rightarrow 
\nonumber\\
& \Rightarrow &
\frac{1}{\hat{a}}\delta_{P}(p,q) - \frac{\hat{c}}{\hat{a}}
\leq  \max\left\{2A\epsilon_{0}, 2\alpha\epsilon_{B}\right\}
\cdot\delta_{\times}(\phi p, \phi q) + AC  
+ \alpha\beta (t_{2}-t_{1})  
\nonumber\\
& \Rightarrow &
\displaystyle{\frac{1}{\hat{a}\cdot
\max\left\{2A\epsilon_{0}, 2\alpha\epsilon_{B}\right\}}}
\delta_{P}(p,q)
- \nonumber\\
& & -\left[\displaystyle{\frac{\hat{c}}{\hat{a}}}
+ AC + \alpha\beta (t_{2}-t_{1})\right]
\displaystyle{\frac{1}
      {\max\left\{2A\epsilon_{0}, 2\alpha\epsilon_{B}\right\}}}
      \leq \delta_{\times}(\phi p, \phi q)
\end{eqnarray}

In what follows, we will produce the inequality that completes
(\ref{leqn442}) into the searched condition (\ref{lequ437}).

\vspace{0.1in}

Let us denote $l:=\delta_{P}(p,q)$.

\vspace{0.1in}

Define  $(y_{0}, y_{1}, \ldots, y_{l})$  a discrete path in $P$ of
minimum length $l$ joining $p$ to $q$. 
Hence,  $(y_{0}, y_{1},
\ldots, y_{l})$ has the following properties,
\begin{eqnarray}
\protect\label{leqn443}
y_{0}:= p \in P, & & y_{l}:= q \in P,\nonumber\\
 d_{M}(y_{\imath},y_{\jmath})\geq \hat{\epsilon},& &
\imath,\jmath=0,1,\ldots,l \hspace{0.1in}(\imath\neq\jmath)\nonumber\\
d_{M}(y_{\iota-1},y_{\iota})\leq 2\hat{\epsilon}, & &
\iota=1,\ldots,l
\nonumber\\
\Rightarrow \delta_{P}(y_{\iota-1},y_{\iota})=1, & &
\iota=1,\ldots,l
\nonumber\\
\stackrel{\mbox{\scriptsize def }P}{\Rightarrow}\pi y_{\iota}\in P_{B},
& & \iota=0,1,\ldots,l \nonumber\\
\stackrel{\mbox{\scriptsize def }P}{\Rightarrow}
\varphi_{(\gamma_{[\pi y_{\iota},b_{0}]})}y_{\iota}\in P_{0},
& & \iota=0,1,\ldots,l \nonumber\\
\stackrel{\mbox{\scriptsize def } P_{B}}{\Rightarrow}
d_{B}(\pi y_{\imath},\pi y_{\jmath})\geq\epsilon_{B}, & &
\imath,\jmath=0,1,\ldots,l \hspace{0.1in}(\imath\neq\jmath)
\end{eqnarray}

Next, we will compare $l$ with $\delta_{P_{B}}\left(\pi p,\pi q\right)$.

\vspace{0.1in}

Notice that because we are assuming ({\bf HLC}), by  
{\bf Lemma~\ref{llem32}} we obtain for any $x,y \in M$,
\begin{equation}
\protect\label{lequ444}
d_{M}(x,y)
\geq \displaystyle{\frac{1}{\alpha}}d_{B}(\pi x,\pi y)-\beta
\end{equation}

So with (\ref{leqn443}), (\ref{lequ444}) and
(\ref{lequ463})%{\bf Lemma~\ref{llem51}} ({\it 8.})
, we get,
\begin{eqnarray}
\protect\label{leqn445}
&\Rightarrow & 2\hat{\epsilon}\stackrel{(\ref{leqn443})}{\geq} 
d_{M}(y_{\iota-1},y_{\iota}) \stackrel{(\ref{lequ444})}{\geq} 
\displaystyle{\frac{1}{\alpha}}
d_{B}(\pi y_{\iota-1},\pi y_{\iota})-\beta \Rightarrow \nonumber\\
&\Rightarrow & \alpha(2\hat{\epsilon}+\beta) \geq
d_{B}(\pi y_{\iota},\pi y_{\iota-1}) \Rightarrow \nonumber\\ 
&\Rightarrow & 2\epsilon_{B}
\stackrel{%\mbox{\scriptsize \bf Lem~\ref{llem51}}({\it 8.})
(\ref{lequ463})}{>}
\alpha(2\hat{\epsilon}+\beta) \geq
d_{B}(\pi y_{\iota},\pi y_{\iota-1})
\stackrel{(\ref{leqn443})}{\geq} \epsilon_{B} \Rightarrow \nonumber\\ 
&\Rightarrow & 2\epsilon_{B} \geq
d_{B}(\pi y_{\iota},\pi y_{\iota-1}) \geq \epsilon_{B}, 
\hspace{0.4in} \forall \iota=1,\ldots,l 
\end{eqnarray}

Since (\ref{leqn445}) holds, we obtain a discrete path
$(\pi y_{0},\pi y_{1},\ldots,\pi y_{l-1},\pi y_{l})$ in $P_{B}$,
connecting $\pi p$ to $\pi q$. 

Therefore, by the definition of $\delta_{P_{B}}$, we conclude that,
\begin{equation}
\protect\label{lequ446}
\delta_{P_{B}}(\pi p, \pi q)\leq l=\delta_{P}(p,q)
\end{equation}

Now, we will compare $l$ with $\delta_{P_{0}}
\left(\varphi_{(\gamma_{[\pi p,b_{0}]})}p,
\varphi_{(\gamma_{[\pi q,b_{0}]})}q \right)$.

By Lemma~2.5~\cite{MK1} , we have for the nets $P_{0}$, $P$ and 
$P_{B}$, respectively,
\begin{eqnarray}
\protect\label{leqn447}
\lefteqn{\exists a_{0}(m,k_{M},\epsilon_{0}) \geq 1,
\exists c_{0}(m,k_{M},\epsilon_{0}) > 0:} \nonumber\\
& & \delta_{P_{0}}\left(p_{1}, p_{2}\right)  
\leq a_{0}\cdot d_{M}\left(p_{1}, p_{2}\right) + c_{0},\hspace{0.2in}
\forall p_{1}, p_{2}\in P_{0}
\end{eqnarray}
\begin{equation}
\protect\label{lequ448}
d_{M}\left(p_{3}, p_{4}\right)\leq
2\hat{\epsilon}\delta_{P}\left(p_{3}, p_{4}\right),\hspace{0.2in}
\forall p_{3}, p_{4}\in P
\end{equation}
\begin{equation}
\protect\label{lequ449}
d_{B}(\pi y_{\imath-1}, \pi y_{\imath}) \leq 2\epsilon_{B}
\delta_{P_{B}}(\pi y_{\imath-1}, \pi y_{\imath}),\hspace{0.2in}
\forall \imath=1,\ldots,l 
\end{equation}

Initially, for each $\imath=1,\ldots,l$, let us look at
\[
\delta_{P_{0}}
\left(\varphi_{(\gamma_{[\pi y_{\imath-1},b_{0}]})}y_{\imath-1},
\varphi_{(\gamma_{[\pi y_{\imath},b_{0}]})}y_{\imath}\right)
\]
which is well defined because of properties 
(\ref{leqn443}).

For each unique (see remark in the beginning of {\it Step 1.})
minimal geodesic in $B$ joining $\pi y_{\imath-1}$ to $\pi y_{\imath}$,
let its unique horizontal lift through $y_{\imath-1}$ be denoted by
$\Gamma_{y_{\imath-1}}: [t_{1}, t_{2}]\longrightarrow M$.

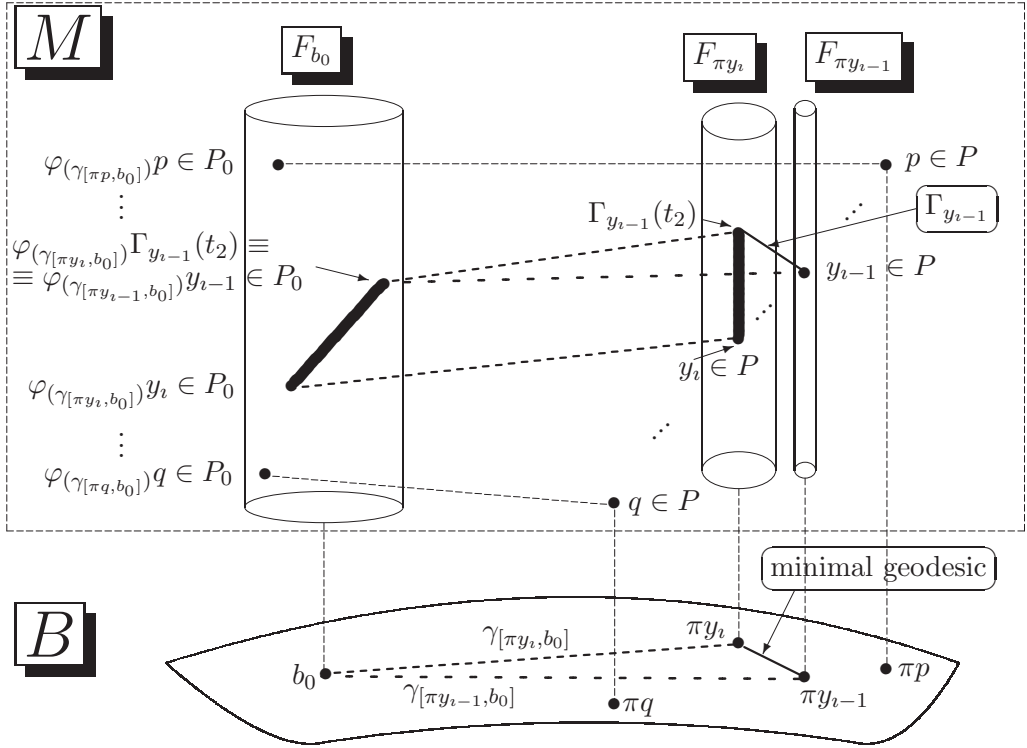
\begin{figure}[here]

         \begin{picture}(390,286)(0,0)%(380,260)(0,10)%

%\dottedline{2}(0,286)(380,286)
%\dottedline{2}(0,0)(380,0)

\put(3,256){\shadowbox{\Huge $M$}}
\dashline[+90]{3}(0,285)(385,285)
\dashline[+90]{3}(0,85)(0,285)
\dashline[+90]{3}(385,85)(385,285)
\dashline[+90]{3}(0,85)(385,85)

\put(104,254){\shadowbox{\large $F_{b_{0}}$}}
%cylinder
%parallels=horizontal curves
\put(120,244){\ellipse{60}{12}}
\put(120,94){\ellipse{60}{12}}
%meridians=vertical curves
\put(90,94){\line(0,1){150}}
\put(150,94){\line(0,1){150}}

%\thicklines
%\put(100,110){\line(-1,3){38}}
%\thinlines

%net P_{0}
\put(100,221){$\bullet$}
\put(14,222){$\varphi_{(\gamma_{[\pi p,b_{0}]})}p\in P_{0}$}
\dashline[+90]{3}(102,224)(330,224)%dashline \varphi_{(\gamma_{[\pi p,b_{0}]})}p --p
\put(42,203){$\vdots$}%\put(58,200){$q_{1}$}
\put(140,176){$\bullet$}%varphi yi-1
\put(2,190){$\varphi_{(\gamma_{[\pi y_{\imath},b_{0}]})}
             \Gamma_{y_{\imath-1}}(t_{2})\equiv$}
\put(2,178){$\equiv\varphi_{(\gamma_{[\pi y_{\imath-1},b_{0}]})}
             y_{\imath-1}\in P_{0}$}%varphi yi-1
\put(117,185){\vector(4,-1){20}}
\thicklines
\dashline[+90]{3}(141,179)(274,198)%dashline varphi yi Gamma yi-1(t2)--Gamma yi-1(t2) 
\thinlines
\thicklines
\dashline{3}[0.7](141,179)(302,183)%dashline varphi yi-1--yi-1
\thinlines
%%\put(124,210){$\bullet$}
%%\dashline[+90]{3}(126,213)(280,213)
\dottedline[$\bullet$]{1}(109,141)(142,178)%varphi yi-1--varphi yi
\put(105,137){$\bullet$}%varphi yi
\put(8,137){$\varphi_{(\gamma_{[\pi y_{\imath},b_{0}]})}
              y_{\imath}\in P_{0}$}%varphi yi
\thicklines
\dashline[+90]{3}(107,139)(275,158)%dashline varphi yi--yi
\thinlines
\put(42,113){$\vdots$}
\put(95,104){$\bullet$}
\put(14,104){$\varphi_{(\gamma_{[\pi q,b_{0}]})}q\in P_{0}$}
\dashline[+90]{3}(97,106)(227,95)%dashline \varphi_{(\gamma_{[\pi q,b_{0}]})}q\in P_{0}--q 

%\put(200,270){$\varphi_{(\gamma_{[\pi p,b_{0}]})}$}
%\thicklines
%\put(302,260){\vector(-1,0){167}}
%\thinlines

\put(256,251){\shadowbox{\large $F_{\pi y_{\imath}}$}}
%cylinder
%parallels=horizontal curves
\put(277,240){\ellipse{28}{14}}
\put(277,108){\ellipse{28}{14}}
%meridians=vertical curves
\put(263,108){\line(0,1){132}}
\put(291,108){\line(0,1){132}}

\dashline[+90]{3}(277,41)(277,101)

\put(302,251){\shadowbox{\large $F_{\pi y_{\imath-1}}$}}
%cylinder
%parallels=horizontal curves
\put(302,245){\ellipse{8}{6}}
\put(302,108){\ellipse{8}{6}}
%meridians=vertical curves
\put(298,108){\line(0,1){137}}
\put(306,108){\line(0,1){137}}

\dashline[+90]{3}(302,28)(302,105)

%net P
\put(330,221){$\bullet$}
\put(340,223){$p\in P$}
\thicklines
\dottedline[$\cdot$]{4}(318,203)(323,208)%p--yi-1
\thinlines
%\put(340,203){$\vdots$}
\put(274,195){$\bullet$}%\Gamma_{y_{\imath-1}}(t_{2})
\put(220,202){$\Gamma_{y_{\imath-1}}(t_{2})$}
\put(263,203){\vector(4,-1){11}}
\dottedline[$\bullet$]{1}(277,197)(277,159)
\put(344,205){\ovalbox{$\Gamma_{y_{\imath-1}}$}}
\put(344,207){\vector(-4,-1){56}}
\thicklines
\put(278,199){\line(3,-2){24}}
\thinlines
\put(299,180){$\bullet$}
%\put(337,175){\vector(-4,1){31}}
%\put(340,174){$y_{\imath-1}\in P$}
\put(310,182){$y_{\imath-1}\in P$}
\thicklines
\dottedline[$\cdot$]{4}(284,165)(289,170)%yi-1--yi
\thinlines
\put(274,155){$\bullet$}
%\put(337,144){\vector(-4,1){55}}
%\put(340,143){$y_{\imath}\in P$}
\put(260,150){\vector(3,1){15}}
\put(254,144){$y_{\imath}\in P$}
\thicklines
\dottedline[$\cdot$]{4}(245,120)(250,125)%yi--q
\thinlines
\put(227,93){$\bullet$}
%\put(337,95){\vector(-1,0){102}}
\put(235,93){$q\in P$}
%\put(340,104){$\vdots$}
%\put(335,107){$\varphi^{-1}_{(\gamma_{\pi p})}q_{l}$}
%\put(343,90){$\in P$}

\put(3,30){\shadowbox{\Huge $B$}}
\qbezier(60,35)(220,85)(360,35)%\qbezier(170,25)(247,80)(340,25)
\qbezier(60,35)(100,0)(120,5)%\qbezier(170,25)(195,25)(210,5)
\qbezier(120,5)(220,20)(320,5)%\qbezier(210,5)(247,20)(300,5)
\qbezier(360,35)(340,0)(320,5)%\qbezier(340,25)(315,25)(300,5)

\put(118,28){$\bullet$}
\put(108,27){$b_{0}$}
\dashline[+90]{3}(120,32)(120,88)
\put(227,17){$\bullet$}
\put(233,17){$\pi q$}
\dashline[+90]{3}(230,19)(230,93)%\dashline[+90]{3}(196,29)(196,109)
\put(330,30){$\bullet$}%\put(289,37){$\bullet$}
\put(337,30){$\pi p$}
\dashline[+90]{3}(333,31)(333,223)%\dashline[+90]{3}(291,39)(291,223)
\put(274,40){$\bullet$}
\put(257,46){$\pi y_{\imath}$}
\put(180,44){$\gamma_{[\pi y_{\imath},b_{0}]}$}
\thicklines
%%\put(279,43){\line(-13,-1){160}}%pi yi--bo
\dashline[+90]{3}(119,31)(275,42)%dashline pi yi--bo
\thinlines
\put(299,27){$\bullet$}
\put(300,20){$\pi y_{\imath-1}$}
\put(150,22){$\gamma_{[\pi y_{\imath-1},b_{0}]}$}
\thicklines
%%\put(300,30){\line(-1,0){180}}%pi yi-1--bo
\dashline{3}[0.7](119,30)(300,29)%dashline pi yi-1--bo
%\dottedline{2}(119,30)(300,29)
\put(279,41){\line(2,-1){24}}%pi yi--pi yi-1
\thinlines

%text

\put(285.5,68){\ovalbox{minimal geodesic}}
\put(298,62.5){\vector(-1,-2){12}}
   \end{picture}

         \caption{$y_{\imath-1}$, the lift
         $\Gamma_{y_{\imath-1}}(t_{2})$, and their
         counterparts in the fiber $F_{b_{0}}$.}
         \label{fig416.2thm53}
         \index{pictures!Theorem\ref{lthm53}}
      \end{figure}

By trivial holonomy {\bf (TH)}, (\ref{leqn447}) 
and {\bf (RIF)}, we may write (see Fig.~{\ref{fig416.2thm53}})
\begin{eqnarray}
\protect\label{leqn450}
\lefteqn{\delta_{P_{0}} 
\left(\varphi_{(\gamma_{[\pi y_{\imath-1},b_{0}]})}y_{\imath-1},
\varphi_{(\gamma_{[\pi y_{\imath},b_{0}]})}y_{\imath}\right)\leq}
\nonumber \\
& \stackrel{{\bf (TH)}}{\leq} &
\delta_{P_{0}}\left(
\varphi_{(\gamma_{[\pi y_{\imath},b_{0}]})}\Gamma_{y_{\imath-1}}
(t_{2}),
\varphi_{(\gamma_{[\pi y_{\imath},b_{0}]})}y_{\imath}
\right) \leq \nonumber \\  
& \stackrel{(\ref{leqn447})}{\leq} &
a_{0}d_{M}\left(
\varphi_{(\gamma_{[\pi y_{\imath},b_{0}]})}\Gamma_{y_{\imath-1}}
(t_{2}),
\varphi_{(\gamma_{[\pi y_{\imath},b_{0}]})}y_{\imath}
\right)  + c_{0} \leq \nonumber \\ 
& \stackrel{{\bf (RIF)}}{\leq} &
a_{0}\left[
A d_{M}\left(\Gamma_{y_{\imath-1}}(t_{2}),y_{\imath}\right) + C
\right]  + c_{0} 
\end{eqnarray}

%(see Fig.~{\ref{fig416.2thm53}})

Now, that (\ref{lequ448}), {\bf Proposition~\ref{lprop35}} and
(\ref{lequ449}) hold, implies
\begin{eqnarray}
\protect\label{leqn451}
\lefteqn{d_{M}\left(\Gamma_{y_{\imath-1}}(t_{2}), y_{\imath}\right)
\stackrel{\triangle}{\leq}
d_{M}\left(\Gamma_{y_{\imath-1}}(t_{2}), y_{\imath-1}\right) +
d_{M}\left(y_{\imath-1}, y_{\imath}\right) \leq} \nonumber \\
& \stackrel{(\ref{lequ448})}{\leq} &
d_{M}\left(\Gamma_{y_{\imath-1}}(t_{2}), y_{\imath-1}\right) +
2\hat{\epsilon}\overbrace{\delta_{P}(y_{\imath-1}, y_{\imath})}^{
   \stackrel{(\ref{leqn445})}{=}\bf{1}} 
%\leq  \nonumber \\
 \stackrel{dist.}{\leq}  \ell(\Gamma_{y_{\imath-1}}) +
2\hat{\epsilon} \leq  \nonumber \\
& \stackrel{Prop\ref{lprop35}}{\leq} &
\alpha\left[\ell\left(\gamma_{\pi y_{\imath-1}\pi y_{\imath}}\right)
+ \beta(t_{2}-t_{1})\right] +
2\hat{\epsilon} =  \nonumber \\
& \stackrel{dist.}{=} & \alpha d_{B}(\pi y_{\imath-1},\pi y_{\imath})
+ \alpha\beta(t_{2}-t_{1}) +
2\hat{\epsilon}\leq  \nonumber \\
& \stackrel{(\ref{lequ449})}{\leq} &
\alpha[2\epsilon_{B}
\underbrace{\delta_{P_{B}}(\pi y_{\imath-1},\pi y_{\imath})}_{
\stackrel{(\ref{leqn445})}{=}\bf{1}}]
+ \alpha\beta(t_{2}-t_{1}) +
2\hat{\epsilon} =\nonumber\\
& \stackrel{(\ref{leqn445})}{=} & 
\alpha2\epsilon_{B} + \alpha\beta(t_{2}-t_{1}) +
2\hat{\epsilon}
\end{eqnarray}

By combining and
(\ref{leqn450}) and (\ref{leqn451}), we get
\begin{eqnarray}
\protect\label{leqn452}
\lefteqn{\delta_{P_{0}} 
\left(\varphi_{(\gamma_{[\pi y_{\imath-1},b_{0}]})}y_{\imath-1},
\varphi_{(\gamma_{[\pi y_{\imath},b_{0}]})}y_{\imath}\right)\leq} 
\nonumber\\
& \stackrel{(\ref{leqn450})}{\leq} & 
a_{0}\left[
Ad_{M}\left(\Gamma_{y_{\imath-1}}(t_{2}),y_{\imath}\right) + 
C\right] + c_{0} \leq \nonumber\\
&  \stackrel{(\ref{leqn451})}{\leq} &
a_{0}\left\{A\left[
2\alpha\epsilon_{B} + \alpha\beta(t_{2}-t_{1}) + 
2\hat{\epsilon}\right]
+ C\right\} + c_{0} = \nonumber\\
& = & a_{0} A\left[
2\alpha\epsilon_{B} + \alpha\beta(t_{2}-t_{1}) + 
2\hat{\epsilon}\right]
+ a_{0} C + c_{0} 
\end{eqnarray}

If we sum (\ref{leqn452}) over $\imath=1,\ldots,l$, 
we may write
\begin{eqnarray}
\protect\label{leqn453}
\lefteqn{\delta_{P_{0}}
\left(\varphi_{(\gamma_{[\pi p,b_{0}]})}p, 
\varphi_{(\gamma_{[\pi q,b_{0}]})}q \right)\leq}\nonumber\\ 
& \stackrel{\triangle}{\leq} &
\sum_{\imath=1}^{l} \delta_{P_{0}} 
\left(\varphi_{(\gamma_{[\pi y_{\imath-1},b_{0}]})}y_{\imath-1},
\varphi_{(\gamma_{[\pi y_{\imath},b_{0}]})}y_{\imath}\right) \leq
\nonumber \\
&  \stackrel{(\ref{leqn452})}{\leq} &
l\left\{a_{0} A\left[2\alpha\epsilon_{B} + \alpha\beta(t_{2}-t_{1}) + 
2\hat{\epsilon}\right] + a_{0} C + c_{0}\right\} = \nonumber \\
& = & \delta_{P}(p,q)\left\{ 
a_{0} A\left[2\alpha\epsilon_{B} + \alpha\beta(t_{2}-t_{1}) + 
2\hat{\epsilon}\right] + a_{0} C + c_{0}
\right\}
\end{eqnarray}

Next, by using inequalities (\ref{lequ446}) and 
(\ref{leqn453}),
\begin{eqnarray}
\protect\label{leqn454}
\lefteqn{\delta_{\times}(\phi p, \phi q)  =
\delta_{P_{0}}\left(\varphi_{(\gamma_{[\pi p,b_{0}]})}p,
\varphi_{(\gamma_{[\pi q,b_{0}]})}q\right) +
\delta_{P_{B}}(\pi p, \pi q)\stackrel{(\ref{lequ446})}{\leq}} 
\nonumber \\
& \stackrel{(\ref{leqn453})}{\leq} &
\left\{
a_{0} A\left[2\alpha\epsilon_{B} + \alpha\beta(t_{2}-t_{1}) + 
2\hat{\epsilon}\right] + a_{0} C + c_{0}
\right\}\delta_{P}(p,q) 
+ \delta_{P}(p,q) = \nonumber \\
& =& \left\{
a_{0} A\left[2\alpha\epsilon_{B} + \alpha\beta(t_{2}-t_{1}) + 
2\hat{\epsilon}\right] + a_{0} C + c_{0} + 1
\right\} 
\delta_{P}(p,q) 
\end{eqnarray}

Finally, we combine inequalities (\ref{leqn442}) and
(\ref{leqn454}), and  obtain
\begin{eqnarray*}
%\protect\label{leqn437.21}
& & \displaystyle{\frac{1}{\hat{a}\cdot
\max\left\{2A\epsilon_{0}, 2\alpha\epsilon_{B}\right\}}}
\delta_{P}(p,q)
- \nonumber\\
& & -\left[\displaystyle{\frac{\hat{c}}{\hat{a}}}
+ \alpha\beta (t_{2}-t_{1}) + AC\right]
\displaystyle{\frac{1}
      {\max\left\{2A\epsilon_{0}, 2\alpha\epsilon_{B}\right\}}}
\stackrel{(\ref{leqn442})}{\leq} 
\delta_{\times}(\phi p, \phi q) \leq \nonumber \\
& & \stackrel{(\ref{leqn454})}{\leq} 
\left\{a_{0}A\left[2\alpha\epsilon_{B} + 
\alpha\beta(t_{2}-t_{1}) + 2\hat{\epsilon}\right] + a_{0}C + 
c_{0} + 1\right\} 
\delta_{P}(p,q)
\end{eqnarray*}
which is the required property (\ref{lequ437})
\[
\displaystyle{\frac{1}{a}}\delta_{P}(p,q) - c
\leq \delta_{\times}(\phi p, \phi q) \leq a\delta_{P}(p,q) + c
\]
with the universal constants given by
\[
\begin{array}{ll}
a := & \max\left\{
a_{0}A\left[2\alpha\epsilon_{B} + 
\alpha\beta(t_{2}-t_{1}) + 2\hat{\epsilon}\right] + a_{0}C + 
c_{0} + 1,
\right.\\ 
& \left. \hspace{0.45in}\hat{a}\cdot
\max\left\{2A\epsilon_{0}, 2\alpha\epsilon_{B}\right\}\right\}\geq 1\\
c := & \left[\displaystyle{\frac{\hat{c}}{\hat{a}}}
+ \alpha\beta (t_{2}-t_{1}) + AC \right]
\displaystyle{\frac{1}
      {\max\left\{2A\epsilon_{0}, 2\alpha\epsilon_{B}\right\}}} >0
\end{array}
\]
and so $\phi$ satisfies {\bf (RI.2)}.

\vspace{0.1in}

This concludes the proof of the Theorem.

\pfe

%\vspace{0.2in}

\begin{center}
{\bf Acknowledgments}
\end{center}

%I wish to express my gratitude to everyone who contributed 
%to making this article a reality.

I wish to register my sincere gratitude and thanks to
Professor Edgar Feldman, my doctoral thesis advisor, and to 
Professor Christina Sormani for assistance with exposition.

  \nocite{*}

  \bibliography{mxrankbgarxiv}

\begin{flushright}
%C. Abreu-Suzuki\\
Mathematics, CIS Department\\
SUNY College at Old Westbury\\
Old Westbury, NY 11568\\
email: casuzuki@earthlink.net\\
\end{flushright}

\end{document}